\documentclass[11pt,reqno]{amsart}
\usepackage[french]{babel} 
\usepackage[latin1]{inputenc}  
\usepackage{amssymb}
\usepackage{stmaryrd}
\usepackage[all]{xy}
\usepackage{graphicx}
\usepackage{hyperref}
\theoremstyle{plain}

\newtheorem{theore}{Théorème}
 
\newtheorem{theo}{Théorème}[section]
\newtheorem{lemme}[theo]%
 {Lemme}
{Définition}
\newtheorem{prop}[theo]%
{Proposition}
\newtheorem{coro}[theo]%
{Corollaire}
{Définition-Proposition}
{Conjecture}

\newcommand{\finpreuve}{\mbox{} \hfill \mbox{$\Box$}}
\newcommand{\finpreuvelemme}{\mbox{} \hfill \mbox{$\blacksquare$}}
\newenvironment{preuve}{\noindent {\it Preuve:}}{\finpreuve}
\newenvironment{preuvelemma}{\noindent {\it Preuve du lemme:
}}{\finpreuvelemme}

\makeatletter
\newcounter{enum@ux}
\def\brkenum#1{%
\setcounter{enum@ux}{\value{enum\romannumeral\the\@enumdepth}}%
\end{enumerate} #1%
\begin{enumerate}%
\addtocounter{enum\romannumeral\the\@enumdepth}{\value{enum@ux}}}
\makeatother

\newcommand{\pg}{\mathfrak{p}}
\newcommand{\Pg}{\mathfrak{P}}
\newcommand{\qg}{\mathfrak{q}}
\newcommand{\Qg}{\mathfrak{Q}}
\newcommand{\tg}{\mathfrak{t}}

\newcommand{\Q}{\mathbb{Q}}
\newcommand{\R}{\mathbb{R}}
\newcommand{\C}{\mathbb{C}}
\newcommand{\QQ}{\mathcal{Q}}
\newcommand{\F}{\mathbb{F}}
\newcommand{\FF}{\mathbb{F}_r(t)}
\newcommand{\KK}{\mathcal{K}}

\newcommand{\N}{\mathrm{N}}

\title[Invariants de Tsfasman-Vl\u adu\c t]{Quelques résultats effectifs concernant les invariants de Tsfasman-Vl\u adu\c t}
\author{Philippe Lebacque}
\date{}
\address{School of Mathematical Sciences, University of Nottingham, Nottingham
 NG7 2RD,
United Kingdom}
\email{Philippe.Lebacque@nottingham.ac.uk}
\begin{document}
\maketitle
\begin{abstract} On considère dans cet article les propriétés asymptotiques de corps globaux à travers l'étude de leurs invariants de Tsfasman-Vl\u adu\c t, nombres qui décrivent en particulier la décomposition des places dans les tours de corps globaux. On utilise des résultats récents de Schmidt et une version faible mais effective du théorème de Grunwald-Wang pour construire des corps globaux infinis ayant un ensemble fini donné d'invariants non nuls et un ensemble prescrit d'invariants nuls, tout en estimant leur défaut.
\end{abstract}










\renewcommand{\thetheo}{\Alph{theo}}

Dans les années 1980, Ihara (voir \cite{IPD}) a initié la théorie asymptotique des corps de nombres, en s'interrogeant sur le nombre de places pouvant se décomposer dans une extension algébrique infinie non ramifiée d'un corps de nombres, et précisa alors très fortement le théorème de densité de Cebotarev, qui prévoit que ces places ont une densité analytique nulle. Ce problème est en outre très important dans le cas des corps de fonctions, puisqu'il est lié à la recherche des courbes ayant un très grand nombre de points rationnels, courbes utiles à la théorie des codes ou encore dans les problèmes d'empilement de sphères, où l'on s'intéresse à la construction de familles de corps de fonctions dont la limite du nombre de points rationnels sur le genre (plus précisément du ratio $N_q/(g-1)$) est maximale (voir \cite{TV}). 
Drinfeld et Vl\u adu\c t (voir \cite{DV}) ont démontré que pour toute famille de courbes sur $\F_q,$ la limite supérieure du ratio $N_q/g$ ne pouvait excéder $\sqrt{q}-1,$ améliorant ainsi la borne obtenue directement par l'application de l'inégalité de Hasse-Weil. Différentes approches (voir \cite{GS}) permettent de construire des familles de courbes sur $\F_{q^2}$ atteignant cette borne, ou d'obtenir des familles sur $\F_q$ dont cette limite est positive.

Tsfasman et Vl\u adu\c t ont par la suite généralisé la borne de Drinfeld-Vl\u adu\c t et les travaux d'Ihara aux familles infinies de corps globaux (et donc aux corps globaux infinis). Leur étude a conduit à des applications diverses, par exemple à une généralisation du théorème de Brauer-Siegel. Ils ont ainsi défini un ensemble d'invariants dont l'importance se voit dans la fonction zêta des corps globaux infinis qu'ils considèrent. Dans sa thèse, l'auteur a tenté de contrôler le support de cet ensemble d'invariants. Si la théorie du corps de classes permet de s'assurer qu'un nombre fini de ces invariants sont positifs, il est plus difficile de répondre au problème inverse: peut-on s'assurer que ces invariants sont nuls. Cela est toutefois rendu possible par les travaux de Labute (voir \cite{LMG}) sur les mild pro-$p$-groupes. Récemment, Schmidt (\cite{SGST}) a généralisé ces résultats aux extensions maximales $S$-ramifiées, $T$-décomposées, et conduit l'auteur à améliorer ses travaux, en contrôlant simultanément un ensemble d'invariants nuls, et un ensemble d'invariants non nuls. C'est ce que nous présentons dans cet article.

\section{Invariants de Tsfasman-Vl\u adu\c t}

On rappelle dans ces paragraphes les définitions et quelques résultats concernant les invariants de Tsfasman-Vl\u adu\c t. On pourra se reporter à \cite{TVF} ou encore à \cite{Lth} pour les détails de ce qui va suivre. 

\subsection{Notations}

Dans toute la suite, on utilisera les conventions et notations suivantes. Par corps global $K$ on entendra une extension finie séparable de $\Q$ ou $\FF,$  pour une puissance $r$ d'un nombre premier $p.$ On supposera, sauf mention du contraire, que le corps des constantes des corps de fonctions est $\F_r.$  On ajoutera $(CN)$ (respectivement $(CF)$) pour signifier qu'une assertion concerne le cas des corps de nombres (resp. des corps de fonctions).  Dans toute la suite, on désignera par:
\begin{center}
\begin{tabular}{ll}
$\Omega(n)$ & $=\sum \alpha_p$ si $n=\prod p^{\alpha_p}$ est la décomposition de $n$ en facteurs premiers.\\
$\mathcal{Q}$ & le corps $\Q$ (CN), $\FF$ (CF),\\
$\delta_\Q$ & $=1$ (CN), $0$ (CF),\\
$\delta_\F $ & $=1-\delta_\Q$,\\ 
$n_K$ & le degré de $K/\QQ,$\\
$d_K$ & le discriminant de $K$ (CN),\\
$g_K$ & le genre de $K$ $(CF),$ $\log{\sqrt{|d_K|}}$  $(CN)$, appelé également genre de $K,$\\ 
$g^\ast_K$ & $g_K$ $(CN),$ $g_K-1$ $(CF)$\\
$Pl(K)$ & l'ensemble des places de $K,$\\
$Pl_f(K)$ & celui de ses places non archimédiennes,\\
$Pl_r(K)$ & celui de ses places réelles,\\
$\N\pg$ & la norme d'une place $\pg\in Pl_f(K):$ le cardinal du corps résiduel en $\pg,$\\
$\deg\pg$ & $\log\N\pg$ $(CF),$\\
$\Phi_{q}(K)$ & le nombre de places de $K$ de norme $q,$\\
$\Phi_\mathbb{R}(K)$ & le nombre de places réelles de $K,$ \\
$\Phi_\mathbb{C}(K)$ & le nombre de places complexes de $K,$ \\
$\delta_\ell(K)$ & $=1$ si le groupe des racines $\ell^{\text{ème}}$ de $1$ $\mu_\ell\subset K,$ $0$ sinon,\\
$U_v$ & le groupe des unités de $\mathcal{O}_{K_v},$ pour $v\in Pl(K),$\\
& avec la convention $U_v=\mathbb{R}_+^\times $ pour $v$ réelle, $\mathbb{C}^\times$ pour $v$ complexe,\\
$Cl^T(K)$& le groupe des $T$-classes d'idéaux de $K,$\\
$\delta_T(\ell,K)$ & $=1$ si $_\ell Cl^T(K)\neq 0$ (CF), $0$ sinon,\\
$E_{K,T}$  & le groupe des $T$-unités de $K,$\\
$V_S^T(K,\ell)$ & le groupe de Kummer: \\
& $\{a\in K^\times\ |\ a\in K^{\times\ell}_v \text{ pour } v\in S \text{ et } a\in U_vK_v^{\times \ell} \text{ pour } v\notin T\}/K^{\times\ell},$ \\
$K_S^T(\ell)$ & désigne la $\ell$-extension maximale de $K$ non ramifiée hors de $S$ \\ 
& où les places de T sont totalement décomposées,\\ 
$G_S^T(K,\ell)$ & $=Gal(K_S^T(\ell)|K),$\\
\end{tabular} 

\begin{tabular}{ll}
$a_S$& $=\mathrm{pgcd}(\deg\pg,\ \pg\in S)$ (CF), $1$ (CN),\\
$\pi(S)$ & $=\sum_{\pg\in S}\log\N\pg$ pour un ensemble de places finies de $K,$\\
$\pi'(S)$ & $=\log^{+}\pi(S).$
\end{tabular} 
\end{center}
 pour $\ell$ un nombre premier et $S,$ $T\subset Pl(K).$ $\log$ désigne la fonction logarithme en base $e$ dans le cas des corps de nombres, en base $r$ dans le cas des corps de fonctions, et $\log^{+}{x}=\log{x}$ si $x\geq 1$ et $0$ sinon.  On omettra $\ell$ et $K$ dans la notation dès lors qu'aucune ambiguïté n'est à craindre.
 
 Pour une extension $L/K,$ $\pg\in Pl(K)$ et $\Pg\in Pl(L)$ prolongeant $\pg$ à $L,$ on note: 
 
  \begin{center}
\begin{tabular}{ll}
$\Phi_{\pg,q}(L)$ & le nombre de places $\Pg\in Pl(L)$ prolongeant $\pg$ de norme $q,$\\
$S_L\subset Pl(L)$ & l'ensemble des prolongements de places de $S\subset Pl(K)$ à $L,$\\
$S_K\subset Pl(K)$ & l'ensemble des restrictions des places de $S\subset Pl(L)$ à $K.$\\
 \end{tabular} 
\end{center}

\textit{Important:} Enfin, dans ce qui suit, on appellera \textit{constante effective} une constante absolue (ne dépendant d'aucun paramètre) et qu'on peut calculer. On écrira $P\ll Q$ (resp. $P\gg Q$) s'il existe une constante effective $c>0$ telle que $P\leq c\,Q$ (resp. $P\geq c\,Q$). Dans tout ce qui suit, nous mettrons l'accent sur la dépendance des majorations en les différents paramètres (genre, norme des places considérées...), la question de l'optimisation de telles constantes restera alors ouverte après notre étude. 

\subsection{Définition et propriétés des invariants de corps globaux infinis}

Rappelons à présent les définitions relatives aux propriétés asymptotiques des corps globaux. On appelle corps global infini toute extension algébrique infinie de $\QQ$ sans extension des constantes dans le cas des corps de fonctions. Par famille $(K_i)_{i\in\mathbb{N}}$ de corps globaux on désignera une suite $(K_i)_{i\in\mathbb{N}}$ de corps globaux, extensions finies du m\^eme corps de base $(\mathbb{Q}$ ou $\mathbb{F}_r(t)),$ telle que $K_i$ n'est pas isomorphe à $K_j$ si $i\neq j,$ et telle que,  dans le cas des corps de fonctions, le corps des constantes de chaque $K_i$ est le même corps fini $\mathbb{F}_r$ pour tout entier $i.$  Dans une famille, la suite des genres $(g_{K_i})_{i\in\mathbb{N}}$ tend vers l'infini, puisqu'il n'y a qu'un nombre fini de tels corps globaux, à isomorphisme près, de genre plus petit qu'un genre donné $g_0.$ Une tour $\{K_i\}_{i\in\mathbb{N}}$ désignera une famille telle que $K_{i}\subsetneq K_{i+1}$ pour tout $i.$ 

Consid\'erons l'ensemble $$A=\begin{cases}\left\{\mathbb{R},\mathbb{C}, p^k, p \text{ premier, } k\in\mathbb{N}^\ast\right\} \ \qquad (CN)\\ \left\{r^k, k\in \mathbb{N}^\ast\right\}\ \qquad (CF)\end{cases}$$
$A_f$ d\'esignera $A-\{\mathbb{R},\mathbb{C}\}$ dans le cas des corps de nombres, et $A$ dans celui des corps de fonctions.

Soit $\mathcal{K}=\{K_i\}_{i\in \mathbb{N}}$ une famille de corps globaux. 
 On dira que $\mathcal{K}$ est asymptotiquement exacte si, pour tout $q\in A$, la suite $\Phi_q(K_i)/g_{K_i}$ admet une limite, que l'on notera alors $\phi_q(\mathcal{K})$. Dans ce cas, on consid\'erera $\Phi_\mathcal{K}=\left\{\phi_q, q\in A\right\}$.  On omettra $\KK$ dans la notation dès que cela ne prête pas à confusion.
On dira que la famille $\mathcal{K}$ est asymptotiquement bonne si elle est asymptotiquement exacte et qu'au moins l'un des $\phi_q$, $q\in A$, est non nul. Dans le cas contraire on la dira asymptotiquement mauvaise. La limite du ratio $n_{K_i}/g_{K_i}$ si elle existe, sera notée $\phi_\infty(\mathcal{K}),$ à ne pas confondre avec les invariants relatifs à la place $\infty$ dans le cas des corps de fonctions.

Remarquons que de toute famille on peut extraire une famille asymptotiquement exacte. Dans la suite, on ne considérera que des tours $\{K_i\}_{i\in \mathbb{N}}$ de corps globaux, qui sont toujours asymptotiquement exactes (voir \cite{TVF}). De plus, dans ce cas, les limites $\phi_q$ ne dépendent que de la limite $\KK=\cup K_i.$ 
On peut alors définir les invariants de Tsfasman-Vl\u adu\c t $\phi_q(\KK)$ d'un corps global infini $\KK,$ comme étant les $\phi_q$ correspondant à toute tour $\{K_i\}$ telle que $\KK=\cup K_i.$ 

On voit facilement qu'une condition nécessaire à un corps global infini pour être asymptotiquement bon est $\phi_\infty>0,$ cette condition étant également suffisante dans le cas des corps de nombres. Elle est en particulier vérifiée si le corps est non ramifié hors d'un ensemble fini de places, et modérément ramifié sur un corps global (voir \cite{Lth} pour les détails).

Ces invariants vérifient une inégalité fondamentale généralisant l'inégalité de Drinfeld-Vl\u adu\c t (voir \cite{TVF}):
\begin{theore}[Inégalités fondamentales de Tsfasman-Vl\u adu\c t]\label{ineg}  Pour tout corps global infini, on a:
 \begin{align*}
(CN-GRH)\text{  }& \sum_q\frac{\phi_q\log{q}}{\sqrt{q}-1} + (\log{\sqrt{8\pi}}+\frac{\pi}{4}+\frac{\gamma}{2})\phi_\mathbb{R}+(\log{8\pi}+\gamma)\phi_\mathbb{C}\leq 1,\\
(CN)\text{  }\ & \sum_q\frac{\phi_q\log{q}}{q-1}+(\log{2\sqrt{\pi}}+\frac{\gamma}{2})\phi_\mathbb{R}+(\log{2\pi}+\gamma)\phi_\mathbb{C}\leq 1,\\
(CF)\text{  }\ & \sum_{m=1}^\infty \frac{m\phi_{r^m}}{r^{\frac{m}{2}}-1}\leq 1,
\end{align*}
où $\gamma$ est la constante d'Euler, et où l'indication (GRH) signifie une fois pour toute que l'assertion est vraie en supposant l'hypothèse de Riemann généralisée.
 \end{theore}

Pour des raisons pratiques de correspondance entre corps de fonctions et corps de nombres, définissons également, pour toute place $\pg$ d'un corps global $K,$ tout $q\in A$ et tout corps global infini $\KK/K,$ les invariants $\phi_{\pg,q}(\mathcal{K})=\lim \Phi_{\pg,q}({K_i})/g_{K_i},$ pour toute tour $\{K_i\}$ d'extensions de $K$ de réunion $\KK.$ Ces nombres existent et ne dépendent pas de la tour choisie. Le support de $\KK/K$ est alors ainsi défini: 
$$Supp(\KK):=\{\pg\in Pl(\QQ)\ | \ \exists q\in A \quad\phi_{\pg,q}\neq 0\}.$$

On peut définir la fonction zêta d'un corps global infini $\KK$ sous la forme suivante (voir \cite{TVF}):
$$\zeta_\KK(s):=\prod_{q\in A_f}(1-q^{-s})^{-\phi_q}$$ ainsi que sa fonction zêta complétée 
$$(CN)\quad \tilde{\zeta}_\KK(s):=e^s2^{-\phi_\mathbb{R}}\pi^{-s\phi_\mathbb{R}/2}(2\pi)^{-s\phi_\mathbb{C}}\Gamma(s/2)^{\phi_\mathbb{R}}\Gamma(s)^{\phi_{\mathbb{C}}}\zeta_\KK(s),$$ 
$$(CF)\quad  \tilde{\zeta_\KK}:=r^s\zeta_\KK.$$
Le produit eulérien définissant ces fonctions converge absolument pour $Re(s)\geq 1$ ($Re(s)\geq 1/2$ sous $GRH$) d'après les inégalités fondamentales. Il définit alors une fonction analytique sur $Re(s)>1$ ($>1/2$ sous $GRH$). De plus, $\zeta_\KK$ est la limite ponctuelle de $(\zeta_{K_i}^{1/g_{K_i}})$ sur le demi-plan $Re(s)>1$ (voir \cite{TVF}). L'étude de ces fonctions zêta se trouve alors intimement liée à celle des invariants.

Plusieurs questions se posent alors naturellement les concernant. On peut par exemple se demander si l'ensemble des invariants de corps globaux infinis a des propriétés topologiques intéressantes pour une topologie naturelle sur les suites réelles. Toutefois cette question, comme d'autres réputées plus faibles, à savoir si le support des corps globaux infinis peut être infini, sont hors de portée actuellement. En effet, si l'on peut essayer de contrôler le comportement d'un nombre fini de places, aucune technique connue de l'auteur ne permet d'en gérer un nombre infini, tout en veillant à ce que la ramification reste finie (et modérée). On peut également s'interroger sur l'existence de corps globaux infinis ayant un défaut nul, c'est à dire dont la différence entre les deux membres de l'inégalité fondamentale est nulle. De tels corps existent sur $\mathbb{F}_{r^2},$ et peuvent même être obtenus de façon récursive (voir \cite{GS}). Dans le cas des corps de nombres, ou des corps de fonctions sur $\mathbb{F}_r,$ on ne sait pas quelles valeurs peuvent être prises par le défaut.
Posons alors 
\begin{align*}
(CN-GRH)\text{  }\   \delta=1 &-\sum_q\frac{\phi_q\log{q}}{\sqrt{q}-1} 
\\ & - (\log{\sqrt{8\pi}}+\frac{\pi}{4}+\frac{\gamma}{2})\phi_\mathbb{R}-(\log{8\pi}+\gamma)\phi_\mathbb{C},\\
(CN)\text{  }\  \delta=1 &-\sum_q\frac{\phi_q\log{q}}{q-1}\\
& -(\log{2\sqrt{\pi}}+\frac{\gamma}{2})\phi_\mathbb{R}-(\log{2\pi}+\gamma)\phi_\mathbb{C},\\
(CF)\text{  }\  \delta=1 &-\sum_{m=1}^\infty \frac{m\phi_{r^m}}{r^{\frac{m}{2}}-1}.
\end{align*}

Dans un travail précédent (voir \cite{Ltvi}), l'auteur a démontré que, pour toute famille finie de paramètres  $q_1,\dots, q_n\in A_f,$ il existe un corps global infini ayant ses invariants $\phi_{q_1},\dots,\phi_{q_n}$ strictement positifs. Il y est aussi démontré que pour tout ensemble fini de nombres premiers $I,$ il existe un corps de nombres infini asymptotiquement bon $\KK$ tel que $I\cap Supp(\KK)=\emptyset.$ Il est également possible de donner des versions effectives de ces résultats, en terme de défaut du corps obtenu. Toutefois, il n'avait pas pu être démontré qu'il existait un corps global infini ayant ses deux propriétés simultanément. Nous nous proposons alors de démontrer:

\begin{theo}\label{princ}
Soient $P=\{\pg_1,\dots,\pg_n\}\subset Pl_f(\QQ),$ et pour tout $i=1\dots n,$ $n_i$ entiers distincts $d_{i,1},\dots,d_{i,n_i}.$ Soit un ensemble fini $I\subset Pl_f(\QQ)$ tel que $I\cap P=\emptyset.$ On pose $N=\mathrm{ppcm}(n_i)_i\,\mathrm{ppcm}(d_{i,j})_{i,j}.$ Alors il existe un corps global infini $\KK$, totalement réel dans le cas des corps de nombres, tel que: 
\begin{enumerate}
\item $I\cap  Supp(\KK)=\emptyset,$
\item Pour tout $i=1\dots n,$ et tout $j=1\dots n_i,$\quad $\phi_{\pg_i,\N\pg_i^{d_{i,j}}}=\frac{\phi_\infty}{n_id_{i,j}}>0.$
\item $\phi_\mathbb{R}=\phi_\infty$ (CN)
\item Il existe deux fonctions $f(P,N)$ et $g(P,I)$ telles que $$\delta(\KK)\leq 1-\frac{h(P,n_i,d_{i,j})}{f(P,N)+g(P,I)}.$$
 Sous \emph{GRH},  
\begin{align*}
h(P,n_i,d_{i,j}) & =\sum_{i=1}^{n}\frac{\log{\N\pg_i}}{n_i}\sum_{j=1}^{n_i}\frac{1}{\N\pg_i^{\frac{d_{i,j}}{2}}-1}+\delta_\Q\left(\log{8\pi}+\frac{\pi}{4}+\frac{\gamma}{2}\right),
\end{align*}
et si $r$ est premier avec $N$ on peut prendre
\begin{align*}
 f(P,N) &\ll a^{\Omega(N)}\left\{(1+|P|)\log{N} +\pi'(P)+\delta_\F N^2\right\}, \end{align*} où $a>1$ est une constante effective, et 
 \begin{align*}
 g(P,I) & \ll  |P|(\pi'(P)+\log^{+}|I|+\log^{+}\pi'(I))+\pi'(I)\\
 &\quad +(|P|^2+|I|+1)(1+\delta_\F \log a_P)
 \end{align*}
 \end{enumerate}
\end{theo}

On verra qu'on peut prendre ce corps global infini sous la forme d'un compositum $\QQ_S^T(\ell)L,$ où $L$ est une extension finie de $\QQ.$ Le défaut étant d'autant plus petit que la somme $\sum_q\frac{\phi_q\log{q}}{\sqrt{q}-1}$ est grande, on voit que pour rendre le défaut le plus petit possible, il faut que $T$ soit aussi grand que possible pour $\sum_{\pg\in T}\frac{\log{\N\pg}}{\sqrt{\N\pg}-1},$ et $S$ aussi petit que possible pour $\sum_{\pg\in S}{\log{\N\pg}}$ (cette somme intervenant dans le calcul du genre). On peut alors se demander combien de places peuvent se décomposer dans ces extensions $\QQ_S^T(\ell).$ Le théorème de densité de Cebotarev implique que l'ensemble de ces places a une densité analytique nulle. On peut toutefois être plus précis:

\begin{prop}\label{propB}
Soient $\ell$ un nombre premier impair, $S,T\subset Pl_f(\QQ)$ deux ensembles finis disjoints de places finies de $\QQ.$ Dans le cas des corps de nombres, on suppose que $\ell\notin S.$ Dans celui des corps de fonctions, on suppose que $\ell$ est premier à $a_T$ et à $r.$ Soit $D$ l'ensemble des places finies de $\QQ$ totalement décomposées dans $\QQ_S^T(\ell).$ Si $\QQ_S^T(\ell)$ est infinie, alors:
\begin{align*}
(CN)\quad \sum_{\pg\in D}\frac{\log\N\pg}{\N\pg-1}&\leq \frac{\pi(S)}{2}-\log{2\sqrt{\pi}}-\frac{\gamma}{2}\\
(CN-GRH)\quad \sum_{\pg\in D}\frac{\log\N\pg}{\sqrt{\N\pg}-1}&\leq \frac{\pi(S)}{2}-\log{\sqrt{8\pi}}-\frac{\pi}{4}-\frac{\gamma}{2}\\
(CF)\quad \sum_{\pg\in D}\frac{\log\N\pg}{\sqrt{\N\pg}-1}& \leq \frac{\pi(S)}{2}
\end{align*}
\end{prop}
Remarquons que cette inégalité s'avère en pratique trop faible du fait de la taille de $S$ pour montrer qu'aucune autre place que celles de $T$ ne se décomposent totalement. De plus, on se rend compte des limites de nos méthodes actuelles. En effet, les méthodes cohomologiques ne permettent pas pour le moment d'imposer le comportement d'une infinité de places à la fois, tandis que les méthodes analytiques ne permettent pas de détecter un ensemble infini de places décomposées lorsque celui-ci est très petit pour $\sum\frac{\log{\N\pg}}{\sqrt{\N\pg}-1}.$ On est donc impuissant devant le problème de savoir s'il existe ou non d'autres places décomposées dans les extensions $\QQ_S^T(\ell)L.$

L'article s'organise ainsi. Dans la seconde partie, nous construisons l'extension $L,$ au moyen d'un résultat effectif de construction d'extension ayant un comportement local donné, et du théorème de densité de Cebotarev. La troisième partie est consacrée à la construction de l'extension $\QQ_S^T(\ell)$ en utilisant les travaux de Schmidt pour lesquels on donne des versions effectives; on y prouve aussi la proposition \ref{propB}, qui s'avère être un corollaire direct des inégalités fondamentales. Dans la quatrième partie, on compose ces deux constructions pour obtenir le corps voulu et une estimation de son défaut. Enfin on donne en cinquième paragraphe deux estimations valables dans les cas particuliers les plus utiles en pratique. 

L'auteur remercie Alexander Schmidt de l'avoir accueilli à Ratisbonne et d'avoir répondu à ses nombreuses questions. Enfin, ce travail a été en partie financé par la dotation EPSRC EP/E049109  "Two dimensional adelic analysis".

\setcounter{theo}{0}
\renewcommand{\thetheo}{\arabic{section}.\arabic{theo}}

\section{Construction d'un corps global ayant des places de norme donnée}

\subsection{Le théorème de densité de Cebotarev}

Dans ce paragraphe, $L/K$ désigne une extension galoisienne de corps globaux, de groupe de Galois $G$. Dans le cas des corps de fonctions, $L$ et $K$ ont pour corps de constantes  $cste(L)=\F_{r^m}$ et $cste(K)=\F_r$ respectivement. On note $\phi$ la substitution de Frobenius: $x\mapsto x^r.$
Soit $\pi(x)$ la fonction de comptage des places finies de $K$ non ramifiées  dans $L$ de norme inférieure ou égale à $x.$ On pose $\Phi(d)=\pi(r^d)-\pi(r^{d-1})$ dans le cas des corps de fonctions.

 Soit $S$ un ensemble fini d'idéaux premiers de $K,$ au dessus d'un ensemble $S_\mathcal{Q}$ de $\mathcal{Q}.$ Pour un idéal premier $\pg$ de $K,$ le symbole d'Artin $\left(\frac{L/K}{\pg}\right)$ désigne la classe de conjugaison des Frobenius correspondant aux idéaux premiers au-dessus de $\pg$ dans $L.$ Soit $C$ une classe de conjugaison de $G.$ $Pl_f^{nr}$ désignera l'ensemble des places finies de $K$ non ramifiées dans $L/K.$ Posons:
 \begin{align*}
 \pi_C(x) & = \#\left\{\pg\in Pl_f^{nr}\ |\ \N\pg\leq x\text{ et }  \left(\frac{L/K}{\pg}\right)=C\right\}\\
 \Phi_C(d)& =\pi_C(r^d)-\pi_C(r^{d-1})\quad (CF)\\
 \end{align*}
  
Considérons la fonction $$\mathrm{Li} (x)=\int_2^x \frac{dt}{\log t}.$$
\begin{theo}[Théorème de densité de Cebotarev \cite{LO},\cite{MS}] Avec les notations précédentes, il existe quatre constantes effectives $A_0,A_1,A_2,A_3$ telles que l'on ait:
$$(CN) \qquad \qquad  \pi_{C}(x)-\frac{|C|}{|G|}\mathrm{Li}(x)= \Delta(x),$$ où, pour tout $x$ tel que $ \log{x}\geq A_0\, n_L\, g_L^2 $
$$(CN)\quad |\Delta(x)|\leq \frac{|C|}{|G|}\mathrm{Li}(x^\rho)+A_1 |\hat{C}|x \exp{\left(-A_2\,n_L^{-\frac{1}{2}}\log^\frac{1}{2}(x)\right)},$$ où $|\hat{C}|$ est le nombre de classes de conjugaisons contenues dans $C,$ et où le terme en $\mathrm{Li}(x^\rho)$ n'est présent que si $\zeta_K$ a un zéro exceptionnel $\rho$ (c'est à dire vérifiant $1-(8\,g_K)^{-1}\leq\rho<1$). 

Sous \emph{GRH}, on a le résultat plus fort suivant, valable pour tout $x\geq 2$:
$$(CN-GRH)\qquad \qquad   |\Delta(x)|\leq A_3\frac{|C|}{|G|} x^\frac{1}{2}\left(2g_L+n_L\log{x}\right).$$
Dans le cas des corps de fonctions, supposons que $C\subset G$ a pour restriction $\phi^d$ à $cste(L).$ Alors on a:
\begin{align*}
(CF)\quad |\Phi_{C}(d)-\frac{|C|}{|G|}\Phi(d)|\leq & 2 g_L \frac{|C|}{|G|}\frac{r^{\frac{d}{2}}}{d}+2(2g_K+1)|C|\frac{r^{\frac{d}{2}}}{d}\\ 
&+ (1+\frac{|C|}{d})\pi(D),
\end{align*}
où $D= Ram(L/K)$ est l'ensemble des places de $K$ ramifiées dans $L.$ Sinon $\pi_C(d)=0.$
\end{theo}

 \begin{coro} Soit $L/K$ une telle extension galoisienne non triviale. Soit $S$ un ensemble de places finies de $K.$ Alors il existe une place finie $\pg$ de $K$ non ramifiée dans $L$ qui n'est pas dans $S,$ dont le Frobenius est dans la classe $C$ tel que, sous l'hypothèse de Riemann généralisée, on ait:
\begin{align*}
(CN- GRH)\quad \log\mathrm{N}\pg &\leq 2A_4\left\{\log(1+|S_\Q|)+\log{g_L}\right\},\\
(CF)\quad \deg\pg& \leq 2max(A_4\left\{\log(1+|S|)+\log([L:K]+g_L)\right\},m)
\end{align*} pour une constante effective $A_4.$
 \end{coro}
Notons que le terme en $[L:K]$ peut être borné par un terme en $g_L$ dès lors que $g_K>1$ ou que $g_K=1$ et $L/K$ ramifiée. Remarquons également que le $m$ dans l'inégalité provient de l'extension des constantes, qui n'altère pas le genre, mais où peu de places occupent une classe donnée.
 
 \begin{preuve} Il s'agit de trouver $x$ tel que $\pi_C(x)\leq |S|+1.$ Traitons le cas des corps de fonctions, celui des corps de nombres se déduisant encore plus directement de \cite{LO} ou \cite{SC}. D'après la formule de Riemann-Hurwitz pour l'extension galoisienne $L/K,$ on a $$\frac{|G|}{m}g_K+\frac{|G|}{4}\pi(D)< g_L+|G|,$$ De plus, d'après l'inégalité de Weil et la formule d'inversion de Moebius,  on a $$\Phi(d)\geq \frac{r^d}{d}-2(2g_K+1)r^{\frac{d}{2}}.$$ Alors, pour que $\Phi_{C}(d)>|S|,$ il suffit de trouver $d$ tel que :
 \begin{align*}
 r^d\geq &2(2g_K+1)r^{\frac{d}{2}}+2g_Lr^{\frac{d}{2}}+2|G|(2g_K+1)r^{\frac{d}{2}}\\
 & +4\left(\frac{d}{|C|}+1\right)(g_L+\frac{|G|}{m})+|S|,
\end{align*} et tel que $\pi_C(d)\neq 0.$ En jouant sur la constante $A_4$ on obtient le résultat sous la forme voulue.
\end{preuve}

Sans $GRH$, on utilisera un résultat de Lagarias, Montgomery et Odlyzko (\cite{LMO}), qui borne le plus petit premier $\pg$ tel que son Frobenius est dans $C$ par $\log\mathrm{N}\pg\leq A_5g_L.$  On peut de même exclure un ensemble fini de places $S.$

\subsection{Un calcul de genre}

Dans ce paragraphe, nous allons rappeler une majoration du genre d'une extension galoisienne $K/k$ de corps globaux, non ramifiée hors d'un ensemble de places $S,$ lemme très classique (voir \cite{SC} pour l'essentiel) qui nous sera utile à de nombreuses reprises.
\begin{lemme} \label{genre}
Soit $K/k$ une extension de corps globaux, non ramifiée hors d'un ensemble fini de places $S.$ Soit $\ell$ un nombre premier. Dans le cas des corps de fonctions, on suppose $K/k$ modérément ramifiée. Dans celui des corps de nombres, on suppose que seules les places au-dessus de $\ell$ dans $k$ peuvent être sauvagement ramifiées. Alors $$g^\ast_K\leq [K:k] \left(g^\ast_k+\frac{1}{2}\sum_{v\in S}\log{\mathrm{N}v}+\frac{\delta_\Q}{2}n_k\log{[K:k]}\right).$$
Si de plus, $K/k$ est galoisienne, 
$$g^\ast_K\leq [K:k] \left(g^\ast_k+\frac{1}{2}\sum_{v\in S}\log{\mathrm{N}v}+\frac{\delta_\Q}{2} n_k\mathrm{v}_\ell([K:k])\log{\ell}\right)$$
\end{lemme}
\begin{preuve} On applique la formule de Riemann-Hurwitz.
\end{preuve}

\subsection{Une version faible mais effective du théorème de Grunwald-Wang}

Le théorème de Grunwald-Wang prédit l'existence d'extensions abéliennes réalisant un nombre fini de conditions locales. Nous nous proposons modestement de démontrer une version faible effective qui nous suffira. Pourtant il paraît raisonnable qu'une version plus générale puisse être obtenue en construisant un bon corps gouverneur pour le problème général. Dans ce paragraphe, nous allons donc démontrer la proposition suivante: 

\begin{prop}\label{effgrun} Soit $k$ un corps global, $T$ (non vide dans le cas des corps de fonctions) et $I$ deux ensembles disjoints de places finies de $k$ et $\ell$ un nombre premier. Alors il existe une extension abélienne (comprenant éventuellement une extension du corps des constantes) modérément ramifiée $K/k$ d'exposant $\ell,$ telle que les places de $T\cup Pl_r(k)$ sont totalement décomposées et celles de $I$ ont pour degré d'inertie $\ell.$ De plus $g_K$ peut être borné explicitement par une fonction de $\ell,$ $\#I,$ $\#T,$ $n_k$ et $g_k.$ 
\end{prop}

Remarquons que la décomposition des places réelles n'est pas obligatoire: en modifiant la preuve de la proposition, on peut également jouer sur leur ramification. Concernant notre résultat, si $k$ est totalement réel, alors $K$ l'est aussi. On verra qu'on construit $K$ comme le compositum d'une extension non ramifiée $T$-décomposée et d'une extension $\{s\}$-modérément ramifiée, $T$-décomposée. Comme on s'intéresse à l'estimation de $g_K/n_K,$ qui est constant dans les extensions non ramifiées, on ne détériore pas le corps en prenant le compositum par une telle extension non ramifiée.

\subsubsection*{a. Premier cas: la caractéristique $car(k)$ de $k$ est différente de $\ell$}

Rappelons tout d'abord un résultat de réflexion que Georges Gras a établi dans le cas des corps de nombres mais dont la preuve dans le cas des corps de fonctions de caractéristique différente de $\ell$ reste tout à fait valable. Nous ne l'écrirons pas ici, puisque cela se résumerait à copier celle de \cite[V.2.4.4]{GRAS}. Cependant, rappelons quelques unes de ses notations. Pour $k$ un corps global, et $T$ un ensemble fini de places de $k,$ on considère le groupe $V^T=V^T_\emptyset(k,\ell).$ Posons alors $\mathcal{K}_T=k(\sqrt[\ell]{1},\sqrt[\ell]{V^T}).$ On dira qu'une place $v$ est modérée si elle est finie et que $car(k)$ est premier à $\mathrm{N}v.$ Pour deux ensembles de places $S$ et $T$ de $k,$ on dira qu'une extension $K/k$ est $T$-décomposée, $S$-totalement ramifiée, si $K/k$ est non ramifié hors de $S,$ et si toutes les places de $S$ (respectivement de $T$) sont totalement ramifiées (resp. totalement décomposées) dans $K/k.$ Nous pouvons à présent énoncer le résultat de Gras:

\begin{prop}[Gras]\label{GR} Soit $\ell$ un nombre premier et soit $k$ un corps global de caractéristique $p\neq\ell.$ Soit $s$ une place modérée  de $k.$ Soit $T$ un ensemble fini de places de $k$ contenant toutes les places réelles de $k.$ Alors il existe une extension cyclique de degré $\ell$ de $k,$ $T$-décomposée et $\{s\}$-totalement ramifiée si et seulement si $s$ est totalement décomposée dans l'extension $\mathcal{K}_T/k.$
\end{prop}

Nous pouvons à présent démontrer la proposition \ref{effgrun}. \\
$Preuve:$ En plus des notations de \ref{effgrun}, posons $T_0=T\cup Pl_r(k).$ Nous allons considérer l'extension $\ell$-élémentaire non ramifiée $K_1$ de $k$ où toutes les places de $T_0$ sont totalement décomposées. Considérons l'ensemble $I_1\subset I$ des places de $I$ totalement décomposées dans $K_1/k.$ 
\begin{lemme}Pour tout $v\in I_1$, on a $V^{T_0\cup \{v\}}\neq V^{T_0},$ et ainsi  $\mathcal{K}_{T_0}\subsetneq\mathcal{K}_{{T_0}\cup\{v\}}$
\end{lemme}
\begin{preuvelemma}
En effet, si $v\in I_1,$ $v$ est totalement décomposée dans $K_1.$ D'après la théorie du corps de classes, l'image de l'idèle $\textbf{v}$ de $k$ correspondant à $v$ appartient alors à 
$$k^{\times}/k^{\times\ell}\prod_{\pg\in {T_0}}k_\pg^\times/k_\pg^{\times\ell}\prod_{\pg\notin {T_0}}U_\pg/U_\pg^\ell$$ 
et s'écrit $\alpha\textbf{u},$ avec $\alpha\in k^\times/k^{\times\ell}$ et $\textbf{u}\in  \prod_{\pg\in {T_0}}k_\pg^\times/k_\pg^{\times\ell}\prod_{\pg\notin {T_0}}U_\pg/U_\pg^\ell.$ Alors $\alpha=\textbf{v}/\textbf{u} \in V^{{T_0}\cup \{v\}}-V^{T_0}.$ Le fait que les corps soient différents est alors clair d'après la théorie de Kummer.
\end{preuvelemma}

Il nous faut à présent trouver une place $s\notin I_1$ de $k$ ne divisant pas $\ell$ telle qu'il existe une extension cyclique de degré $\ell,$ ${T_0}$-décomposée et $\{s\}$-totalement ramifiée, telle qu'il n'existe pas une telle extension ${T_0}\cup\{v\}$-décomposée, $\{s\}$-totalement ramifiée, pour tout $v\in I_1,$ et de savoir estimer $\mathrm{N}s.$ Cela est rendu possible par les résultats de Gras et le théorème de densité de Cebotarev effectif. En effet, d'après la proposition \ref{GR} une telle place $s$ doit vérifier les conditions suivantes:
\begin{enumerate}
\item $Frob_s\in Gal(\mathcal{K}_{T_0}/k)$ est trivial.
\item Pour tout $v\in I_1,$ $Frob_s\in Gal(\mathcal{K}_{{T_0}\cup\{v\}}/k)$ n'est pas trivial.
\end{enumerate}
Considérons le compositum $L/k$ de toutes les extensions $\mathcal{K}_{{T_0}\cup\{v\}},$ $v\in I_1.$ $L/k$ est une extension galoisienne. Dans $Gal(L/k),$ il existe un élément $\sigma$ trivial sur $\mathcal{K}_{T_0}$ et non trivial sur chacun des $\mathcal{K}_{{T_0}\cup\{v\}},$ puisque chacun de ces corps est différent de $\mathcal{K}_{T_0}$ d'après le lemme. En prenant une place $s$ telle que $Frob_s=\sigma$ on obtient alors, d'après la proposition \ref{GR}, l'existence d'une extension $K_2/k$ cyclique de degré $\ell,$ où toutes les places de $I_1$ sont inertes, et toutes les places de ${T_0}$ sont totalement décomposées. Prenons alors comme extension $K$ le compositum $K_1 K_2.$ Cette extension est abélienne, de groupe de Galois d'exposant $\ell.$ Le degré d'inertie des places de $I$ étant l'ordre du Frobenius correspondant, celui-ci est au plus $\ell,$ et vaut donc $\ell,$ puisque chaque place a pour degré d'inertie $\ell$ dans l'une ou l'autre des deux extensions. Enfin $K$ est non ramifiée hors de $s.$ Reste alors à estimer $\mathrm{N}s.$ Il faut pour cela commencer par estimer le degré et le genre de $L,$ puis on pourra invoquer le théorème de densité de Cebotarev pour obtenir une majoration de la norme de la plus petite place $s$ vérifiant les conditions précédentes. 

\subsubsection*{Estimation du degré et du genre de $L$}
$L$ est le compositum des $\mathcal{K}_{{T_0}\cup\{v\}}/\mathcal{K}_{{T_0}},$ pour tous les $v\in I_1.$ Chacune de ces extensions a pour degré $\ell:$ en effet, elles ne sont pas triviales, et l'application $V^{{T_0}\cup\{v\}}\to \mathbb{Z}/\ell\mathbb{Z},$ qui a un élément $x$ associe sa valuation en $v$ est surjective (sinon $V^{{T_0}\cup\{v\}}=V^{T_0}$) et a pour noyau $V^{{T_0}}.$ 

De plus, le degré de $\mathcal{K}_{{T_0}}/k$ est borné par $$[\mathcal{K}_{{T_0}}:k]\leq[k(\mu_\ell):k]\ell^{d_\ell(V^{{T_0}})},$$
qui, d'après les résultats de Shafarevich \cite{Shaf}, est calculé par :
 $$d_\ell(V^{{T_0}})=d_\ell(Cl^{{T_0}})+\Phi_\mathbb{R}(k)+\Phi_\mathbb{C}(k)-1+\#T+\delta_\ell(k).$$
Posons $a_T=\mathrm{pgcd}(\deg\tg,\tg\in T)$ dans le cas des corps de fonctions, $a_T=1$ dans celui des corps de nombres. Dans ce premier cas, $\#Cl^{\tg}=h_k\deg\tg$ pour tout $\tg\in T$ (voir \cite[1.2.5]{XN}), on a alors:
 $$[L:k]\leq a_T h_k(\ell-1)^{1-\delta_\ell(k)} \ell^{\Phi_\mathbb{R}(k)+\Phi_\mathbb{C}(k)+\delta_\ell(k)-1+\#T+\#I_1}.$$ 
 
La théorie de Kummer affirme de plus que les $\mathcal{K}_{{T_0}\cup\{v\}}/k,$ sont non ramifiées hors de ${T_0}\cup\{v\}\cup P_\ell,$ où $P_\ell$ sont les places de $k$ divisant $\ell,$ seules places pouvant être sauvagement ramifiées. On peut ainsi appliquer le lemme \ref{genre}, et obtenir:

\begin{align*}g^\ast_L\leq & [L:k](g^\ast_k+\frac{1}{2}\sum_{v\in T\cup I_1}\log{\mathrm{N}v})\\ 
& +\frac{\delta_\Q}{2}[L:k]n_k(\#I_1+\#T+  \log_\ell{h_k}+\Phi_\mathbb{R}(k)+\Phi_\mathbb{C}(k))\log{\ell}, \\
\end{align*}
On a alors: 
\begin{align*}
g^\ast_L & \leq  \bar{g}_L:= a_T h_k \ell^{\Phi_\mathbb{R}(k)+\Phi_\mathbb{C}(k)+\#T+\#I_1}\\
& \ \ \times  \left(g_k+\frac{1}{2}\pi(T\cup I_1)+\delta_\Q\frac{\log{\ell}}{2} n_k(\#I_1+\#T+  \log_\ell{h_k}+\Phi_\mathbb{R}(k)+\Phi_\mathbb{C}(k))\right). \\ 
\end{align*}

 \subsubsection*{Majoration de $h_k$} Avant de poursuivre, il nous faut voir que $h_k,$ qui intervient dans la majoration du genre de $L,$ peut également être borné par une fonction de $n_k$ et $g_k.$ 
\begin{lemme} Soit $k$ un corps global. Il existe une constante effective $A_6$ telle que $\log h_k\leq A_6 g_k.$ 
Plus précisément, on a:
\begin{align*}
h_k & \leq \bar{h}_k:=\begin{cases}
25\exp{(-0,46\, n_k)} \left(\frac{e\log |d_k|}{4(n_k-1)}\right)^{n_k-1}\sqrt{|d_k|} \text{ si } k\neq \Q\\
1 \text{ si } k=\Q
\end{cases} (CN)\\
h_k & \leq(1+\sqrt{r})^{2g} \quad (CF)
\end{align*}
\end{lemme}
\begin{preuvelemma}
\textit{Cas des corps de nombres.}
On utilise pour cela les majorations de Louboutin (voir \cite{LOU}) concernant le résidu de la fonction zeta en $1$ d'un corps de nombres $k$ différent de $\Q.$ On se place dans le cas $k\neq\Q.$ On a alors $$h_k \frac{R_k}{w_k}\leq \frac{1}{2}\left(\frac{2}{\pi}\right)^{\Phi_\C(k)} \left(\frac{e\log |d_k|}{4(n_k-1)}\right)^{n_k-1}\sqrt{|d_k|}.$$
Utilisant la minoration de $R_k/w_k\geq 0,02\exp{(0,46\,\Phi_\R(k)+0,1\,\Phi_\C(k))}$ due à Zimmert (voir \cite[\S3]{ZIM}), on obtient pour $h_k$ la majoration:
$$h_k\leq \bar{h}_k:=25\exp{(-0,46\, n_k)} \left(\frac{e\log |d_k|}{4(n_k-1)}\right)^{n_k-1}\sqrt{|d_k|}.$$ Pour $k=\Q$ on prend $\bar{h}_\Q:=1.$
Lorsque le degré est très petit devant le discriminant, on peut utiliser une autre minoration pour le régulateur due à Silverman (voir \cite{Sil}) et gagner un facteur $\log |d_k|$ (à $n_k$ fixé). Si on souhaitait produire une construction asymptotique à base de ce résultat, il serait peut-être intéressant de l'introduire. On se contentera ici de la minoration de Zimmert.

\textit{Cas des corps de fonctions.} Dans ce cas, on peut majorer $h_k$ au moyen de l'hypothèse de Riemann pour les corps de fonctions. En effet, si $P_k(x)$ est le numérateur de la fonction zêta de $k,$ $$h_k=P_k(1)=\prod_{i=1}^{g_k}|1-\rho_i|^2,$$ où $|\rho_i|=\sqrt{r}.$ On obtient alors en majorant très brutalement $$h_k\leq (1+\sqrt{r})^{2g}.$$
\end{preuvelemma}
 
 \subsubsection*{Evaluation de $g_K$}
 
Remarquons d'abord que les corps construits sont des extensions d'exposant $\ell$ de $k(\sqrt[\ell]{1}),$ ainsi l'extension du corps des constantes est au pire de degré $\ell(\ell-1).$ Ici, nous allons avoir besoin du théorème de Cebotarev effectif, pour estimer la norme de la plus petite place $s$ ne divisant pas $\ell$ et se trouvant dans la classe de conjugaison souhaitée.  $K$ est le compositum de $K_1$ et $K_2.$ Ces deux extensions sont linéairement indépendantes, puisque $K_2$ est totalement ramifiée en $s,$ et $K_1$ est non ramifiée ( de sorte que $K/K_2$ est non ramifiée). Puisque $K/k$ est modérément ramifiée, on en déduit, dans le cas des corps de nombres, que: 
 \begin{align*}
 g^\ast_K & =[K_1:k]g^\ast_{K_2}\\ 
 & \leq [K_1:k][K_2:k]g^\ast_k+\frac{1}{2}[K_1:k][K_2:k](1-1/\ell)\log\mathrm{N}s\\
 & \leq  a_Th_k \ell(g^\ast_k+a\bar{g}_L).\\
 &\ll a_Th_k\ell \bar{g}_L,
 \end{align*}
où $a$ est une constante effective.

Ainsi, $g_K$ peut être borné par une fonction explicite de $\ell,$ $\#I,$ $\#T,$ $n_k$ et $g_k,$ par la formule:
\begin{align*} g^\ast_K& \ll a_T^2\bar{h}_k^2\ell^{\Phi_\mathbb{R}(k)+\Phi_\mathbb{C}(k)+1+\#T+\#I_1} \\
&\ \ \ \times \left(g^\ast_k+\frac{1}{2}\pi(T\cup I_1)+\delta_\Q\frac{\log{\ell}}{2} n_k(\#I_1+\#T+ \log_\ell{\bar{h}_k}+\Phi_\mathbb{R}(k)+\Phi_\mathbb{C}(k)\right)
\end{align*}
ceci terminant la preuve du résultat pour les corps de nombres. Toutefois ces majorations sont très mauvaises lorsqu'il s'agit de répéter cette construction et de déterminer le genre du corps ainsi obtenu. Pour ne pas écoeurer le lecteur, nous donnerons plutôt des résultats supposant l'hypothèse de Riemann généralisée.
Si l'on y croit  (elle est vraie pour les corps de fonctions), on obtient alors les résultats suivant (pour le genre normalisé):
\begin{align*}
(GRH)\ \
 \frac{g^\ast_K}{n_K} & =\frac{[K_1:k]}{n_K}g^\ast_{K_2}\\ 
 & \leq \frac{1}{n_k}\left(g^\ast_k+\frac{1}{2}(1-1/\ell)\log\mathrm{N}s\right)\\
 & \leq \frac{1}{n_k}\left\{g^\ast_k+max(A_4(\log(\#I_1+1)+\log{([L:k]+\bar{g}_L)}),\ell(\ell-1)\right\}.\\
 \end{align*} De plus on a  $\log h_k\ll g_k$ et $\log a_T\ll \pi'(T).$ En remplaçant de même $g_L$ par un majorant, on obtient alors:
\begin{align}\label{gK}
(GRH)\ \
 \frac{g_K}{n_K} & \leq A_7\frac{g_k}{n_k}+A_8\left\{\frac{\log{\ell}}{n_k}(\#I_1+\# T+\delta_\Q n_k)+\frac{\log{\pi(T\cup I_1)}}{n_k}\right\}+\delta_\F \frac{\ell^2}{n_k},
  \end{align}
où $A_7$ et $A_8$ sont des constantes effectives. Notons que cette expression est obtenue en prenant la somme plutôt que le maximum dans l'estimation du genre.
\hfill$\square$

 Remarquons que la faiblesse de la majoration (\ref{gK}) provient de l'estimation du degré de $L$ qui introduit de nouveaux termes en $g_k/n_k$ par le biais de $\log{h_k}.$
 
 \ \\
 \textit{b. Cas $\ell=p$.}
 Nous allons utiliser un argument présent en \cite[9.2.5]{NCG2}. Dans le cas d'un corps de fonctions en caractéristique $p=\ell,$ le théorème d'approximation forte prouve que l'application $\mathcal{O}_{k,T\cup I\cup\{\qg\}}\to \oplus_{\pg\in T\cup I}k_\pg/\wp k_\pg$ est surjective, pour toute place $\qg$ non contenue dans $T\cup I,$ avec $\wp(x)=x^p-x.$ En effet, rappelons que, si on se donne, pour chaque $\pg$ d'un ensemble fini $S$ de places d'un corps de fonctions $K,$ des éléments $a_\pg\in K_\pg,$ ainsi que des entiers $n_\pg>0$ et une autre place arbitraire $\qg\notin S,$ il existe un élément $a$ de $K$ tel que $ord_\pg(a-a_\pg)\geq n_\pg$ pour tout $\pg\in S,$ $ord_\pg(a)\geq 0$ pour $\pg\notin S\cup\{\qg\},$ et $ord_\qg(a)\geq 2g+\sum_{\pg\in S}n_\pg$ (voir \cite[6.13]{Rosff} pour ce résultat). Il n'est pas difficile de voir que l'uniformisante $\pi_\pg\in \wp(k_\pg),$ ce qui prouve la surjectivité (prendre $n_\pg=1$).
 
 Soit alors $u\in \mathcal{O}_{k,T\cup I\cup\{\qg\}}$ telle que $u\in \wp(k_\pg)$ pour tout $\pg\in T,$ et $u$ une constante non contenue dans $\wp(k_\pg)$ pour $\pg\in I.$ Le théorème d'approximation forte nous assure qu'on peut prendre $v_{\qg}(u)\geq -(2g_k+|T|+|I|).$ 
 
 Alors, si on considère l'extension d'Artin-Schreier $k(y)$ engendrée par une racine de $y$ de $X^p-X-u,$ on obtient une extension cyclique de degré $p,$ où les places de $T$ sont totalement décomposées, et les places de $I$ sont inertes. De plus elle est non ramifiée hors de $T\cup I\cup\{\qg\},$ c'est à dire qu'elle est non ramifiée hors de $\qg.$ 
 
D'après un calcul classique de genre dans les extensions d'Artin-Schreier (voir \cite[III.7.8]{STIAFFC}, on a:
 $$g_{k(y)}\leq p g_k+\frac{p-1}{2}\left(-2+(2\ g_k+|T\cup I|)\deg\qg\right)$$
 
 Enfin, puisque le nombre $\Phi(d)$ de places de degré $d$ de $k$ est minoré par $\Phi(d)\geq r^d/d-2(2g_k+1)r^{\frac{d}{2}},$
  on peut prendre $$\deg(\qg)\leq 4\log(8+8g_k+|T\cup I|).$$
 
 On en déduit donc, en posant $K=k(y),$ que :
  $$g_K\leq p g_k+(p-1)(2\ g_k+|T\cup I|)(\log_r(3+6g_k+|T\cup I|)).$$

On voit alors que le genre normalisé obtenu est en $g_k\log (g_k|T\cup I|),$ le résultat est alors plus mauvais que dans le cas modéré. On supposera par la suite qu'il n'y a pas de ramification sauvage dans ces extensions.

 \subsection{Construction d'un corps ayant des places de normes données}
 
 Ce paragraphe reprend une construction antérieure de l'auteur (voir \cite{Ltvi}). Nous allons construire un corps global ayant certains $\Phi_{\pg,q}>0.$ Nous estimerons alors son genre. Cette construction pourrait être effectuée à partir de n'importe quel corps global $k,$ toutefois nous nous contenterons du cas où l'on construit une extension de $\QQ,$ puisque c'est cela qui nous intéresse en pratique. 
 
\begin{prop} Soient $\pg_1,\dots,\pg_k$ des places finies de $\QQ$ et $P$ leur ensemble. Soient ${d_{i,1}},\dots, {d_{i,n_i}}$ des entiers positifs pour tout $1\leq i\leq k.$ Alors il existe une extension $L/\QQ$ telle que 
  \begin{enumerate}
\item $\Phi_{\pg_i,\N\pg_i^{d_{i,j}}}(L)=\frac{n_L}{n_i\,d_{i,j}}.$
\item Pour tout $\pg\in P_L$  au dessus d'un des $\pg_i,$ il existe $j$ tel que $\mathrm{N}\pg=\N\pg_i^{d_{i,j}}.$
\item Il existe une fonction $f$ de $P$ et $N=\mathrm{ppcm}(n_i)_i\mathrm{ppcm}(d_{i,j})_{i,j}$ telle que $$g_L/n_L\leq f(P,N),$$ 
pouvant être choisie, sous l'hypothèse de Riemann généralisée, lorsque $N$ est premier avec $r,$ telle que:
\begin{align*}
(GRH)\quad f(P,N) & \ll A_7^{\Omega(N)}\left\{(1+|P|)\log{N} +\pi'(P)+\delta_\F  N^2\right\}.
  \end{align*}
   \end{enumerate}
 \end{prop}
 Nous allons faire la preuve de ce résultat sous \emph{GRH}, la majoration sans \emph{GRH} qu'on obtiendrait étant bien plus faible.
 
 \begin{preuve}
Dans le cas des corps de fonctions, on suppose ici que $r$ est premier à $N.$ De plus, on va ajouter de même que dans la proposition une place dont le degré est premier à $N$ pour s'assurer qu'il n'y a pas d'extension des constantes. Pour cela, il suffit de prendre une place $\qg$ de degré premier $s$ ne divisant pas $N.$ D'après le théorème des nombres premiers, on peut prendre $s\ll \log{N}$ (d'après \cite{HW},on a $\omega(N)\ll \frac{\ln N}{\ln\ln N},$ où $\omega(N)$ est le nombre de facteurs premiers de la décomposition de $N$ comptés sans multiplicité). Dans le cas des corps de nombres, tout cela n'est pas nécessaire. On notera $P'=P\cup \{\qg\}$ dans le cas des corps de fonctions, $P'=P$ pour les corps de nombres.
  
   
Soient $p_1<\dots<p_n$ les nombres premiers divisant l'un des $d_{i,j}.$ Soit $N_1=\mathrm{ppcm}(d_{i,j})_{i,j}.$ 

  On commence alors par considérer une extension $L_0$ de degré $N_2=\mathrm{ppcm}(n_i)$ telle que toutes les places de $P'$ sont totalement décomposées. Pour cela, on applique $\Omega(N_1)$ fois la proposition. Dans $L_0,$ on a ainsi $N_2$ places au-dessus de chaque $\pg_i\in P.$ Pour chacune des places $\pg_i$ de $P,$ on forme $n_i$ ensembles $P_{\pg_i,1},\dots,P_{\pg_i,n_i}$ de $N_2/n_i$ places chacun. Ils sont destinés à fournir de places de norme $\N\pg_i^{d_{i,j}}$ respectivement.
  Puis on construit par récurrence une tour de corps telle que, dans le cas des corps de fonctions, $\qg$ est totalement décomposée (on ne le rappellera pas):
  
 Pour tout $s=1\dots n,$ on construit à partir de $L_{s-1}$ la tour suivante:
 \begin{enumerate}\item $L_{s-1}^0=L_{s-1}$
 \item Pour $t=1\dots v_{p_s}(N_1),$  $L_{s-1}^t/L_{s-1}^{t-1}$ est une extension abélienne d'exposant $p_s$ telle que:\\
pour tout $i,j,$ si $p_s^t | d_{i,j}$ alors les places au-dessus de $P_{\pg_i,j}$ sont inertes. Sinon elles sont totalement décomposées. 
 \item On pose $L_s=L_{s-1}^{v_{p_s}(N_1)}.$
  \end{enumerate}

  Considérons alors $L=L_n,$ et estimons son genre au moyen de (\ref{gK}). Renommons en  $\{K_i\}$ la suite des corps qu'on a construits, $T_i$ l'ensemble des places dont on impose la décomposition, $I_i$ celui des places inertes. $T_i\cup I_i$ est l'ensemble des places au-dessus de $P'$ dans $K_i.$ Ainsi on a 
 $\#I_{1,i}+\#T_i\leq n_{K_i}|P+1|$ à chaque corps intermédiaire $K_i$ construit. De même on a $\log^{+}{\pi(T_i\cup I_i)}= \log{n_{K_i}}+\pi'(P)+\delta_\F \ln\ln{N}.$ 
 On a alors d'après (\ref{gK})
 \begin{align*}
(GRH)\ \
 \frac{g_{K_{i+1}}}{n_{K_{i+1}}} & \leq A_7\frac{g_{K_{i}}}{n_{K_i}}
 +A_8\left\{\frac{\log{\ell_i}}{n_{K_i}}(\#I_{i}+\# T_{i}+\delta_\Q n_{K_i})+\frac{\log^{+}{\pi(T_i\cup I_i)}}{n_{K_i}}\right\}\\
& \quad + \delta_\F \frac{\ell_i^2}{n_{K_i}}\\ 
 & \leq A_7\frac{g^\ast_{K_{i}}}{n_{K_i}}
 +A_8\left\{(1+|P|)\log{\ell_i}+1+\frac{\pi'(P)+\ln\ln{N}}{n_{K_i}}\right\}\\
&\quad +\delta_\F \frac{\ell_i^2}{n_{K_i}},
\end{align*}
 où $\ell_i$ est le nombre premier correspondant à l'extension $K_{i+1}/K_i$ et dépend donc de $i.$
On obtient alors par récurrence immédiate la majoration, si $N=\prod_{i=1\dots m}\ell_i,$ où $\ell_i$ sont les nombres premiers divisant $N:$
\begin{align*}
\emph{(GRH)}\quad \frac{g_L}{n_{L}} & \ll A_7^{\Omega(N)}\left\{(1+|P|)\log{N} +\pi'(P)+\delta_\F \sum_{i=1}^m \frac{\ell_i^2}{\prod_{k<i}\ell_k} \right\}.\\
 \emph{(GRH)}\quad \frac{g_L}{n_{L}} & \ll A_7^{\Omega(N)}\left\{(1+|P|)\log{N} +\pi'(P)+\delta_\F N^2\right\},\\
\end{align*} 
cette dernière inégalité étant réalisée dans le cas le plus mauvais où $N$ est un nombre premier. Comme $\Omega(N)\leq N/\log{2},$ cette inégalité peut encore s'écrire:
 \begin{align*}
\emph{(GRH)}\quad \frac{g_L}{n_{L}} & \ll N^{\frac{A_7}{\log{2}}}\left\{(1+|P|)\log{N} +\pi'(P)+\delta_\F N^2\right\}, \end{align*} 
 \end{preuve}

 On a ainsi achevé la première partie de la construction, c'est à dire qu'on obtient un corps global $L$ dont on sait estimer le genre ayant des $\Phi_{\pg,\N\pg^q}>0$ donnés. Remarquons toutefois que nous avons invoqué un théorème de Gras valable dans le cas d'extensions de degré $\ell.$ Celui-ci admet une généralisation au degré $\ell^i$ que nous n'avons pas utilisé, présumant que cela n'apporterait pas de nette amélioration. Cependant, vu que les majorations les plus coûteuses concernent le nombre de classes, et qu'on se sert de majorations grossières du type $\ell^{d_\ell Cl^T}\leq h_k,$ ils pourraient peut-être améliorer significativement l'estimation du genre (en obtenant une puissance de $N$ inférieure à celle obtenue ici).

 \section{Construction d'un corps global infini dont le support est contrôlé} 
 
 
 \subsection{Résultats quantitatifs relatifs à la propriété $K(\pi,1)$ (voir \cite{SGST})}

 Le but de cette section est de donner une version quantitative de résultats de ce type:
 
 \begin{theo}[Schmidt, voir \cite{SGST}]\label{Schmidt1}
 Soient $T$ et $S$ deux ensembles finis disjoints de places finies d'un corps global $k,$ tel que $T$ soit non vide dans le cas des corps de fonctions. Soit $\ell$ un nombre premier impair différent de $car(k).$ Alors il existe un ensemble fini de places finies $S_0$ ne contenant pas de places divisant $\ell,$ tel que $cd\,G(k_{S\cup S_0}^T|k)(\ell)=2$ et que pour toute place $\pg\in S\cup S_0,$ $k^T_{S\cup S_0}(\ell)_\pg=k_\pg(\ell),$ c'est à dire que l'extension  $k^T_{S\cup S_0}(\ell)$ réalise la $\ell$-extension maximale  $k_\pg(\ell)$ de $k_\pg.$
 \end{theo}

 Nous allons reprendre les grandes lignes de la preuve afin d'obtenir une précision sur la taille des places de $S_0,$ et ainsi nous pourrons obtenir des informations sur les invariants de ce corps. A présent $\ell$ désignera un entier impair premier à la caractéristique de $k.$
 
 \subsubsection{Quelques définitions} Rappelons tout d'abord la notion de courbe marquée et de son site étale comme on peut les trouver dans \cite{SGST}. Pour un schéma noethérien régulier $Y$ de dimension $1$ et $T$ un ensemble fini de points fermés de $Y,$ Et($Y$) désigne à l'accoutumée  la catégorie des morphismes étales de type fini $Y'\to Y,$ et Et($Y,T$) la sous-catégorie pleine dont les objets sont les morphismes $f:Y'\to Y$ tels que, pour tout point fermé $y'\in Y'$ tel que $f(y')=y\in T,$ l'extension résiduelle $k(y')|k(y)$ est triviale. Le site étale $(Y,T)_{et}$ consiste alors en la catégorie Et($Y,T$) munie des familles surjectives comme recouvrements. De même que pour le site étale, on construit le groupe fondamental, que l'on notera $\pi_1^{et}(Y,T,\bar{x})$ (où $\bar{x}$ est le point géométrique de $Y-T$ choisi pour la construction), ou plus simplement $\pi_1^{et}(Y,T)$ lorsque $Y$ est connexe. Ce groupe classifie les recouvrements étales où les points de $T$ sont totalement décomposés. Si on considère le pro-$\ell$-recouvrement universel $\widetilde{(Y,T)}(\ell)$ de $(Y,T),$ on a alors la suite spectrale 
 $$E_2^{ij}=H^i\left(\pi^{et}_1(Y,T)(\ell),H^j_{et}(\widetilde{(Y,T)}(\ell),\F_\ell)\right)\Rightarrow H^{i+j}_{et}(Y,T,\F_\ell),$$ et en particulier les edge morphismes 
 $$\phi_{i}:H^i\left(\pi^{et}_1(Y,T)(\ell),\F_\ell\right)\to H^i_{et}(Y,T,\F_\ell),\quad i\geq 0.$$ 
 Pour $i=0,1$ ce sont alors des isomorphismes, et pour $i=2$ le morphisme est injectif. S'ils sont des isomorphismes pour tout $i$ (ou de façon équivalente, $H^j_{et}(\widetilde{(Y,T)}(\ell),\F_\ell)=0$ pour tout $j\geq 1$) on dira que $(Y,T)$ jouit de la propriété $K(\pi,1)$ pour $\ell.$
 
 Dans le cas qui nous intéresse, $k$ est un corps global, $Y$ est $X-S,$ où $X=Spec\,\mathcal{O}_k,$ et $S$ est un ensemble fini de places finies de $k.$ Alors on a $\pi^{et}_1(Y,T)(\ell)=Gal(k_{S\cup Pl_r}^T(\ell)|k)=Gal(k_{S}^T(\ell)|k)$ car $\ell$ est impair. Si $K|k$ est une sous-extension de $k^T_S(\ell)|k,$ on note $(X-S,T)_K$ la normalisation de la courbe $X-S$ dans $K$ marquée aux points de $T_K.$
 
Si $S\neq\emptyset,$ les assertions suivantes sont alors équivalentes (voir \cite{SGST}):
 \begin{enumerate}
 \item $(X-S,T)_K$ a la propriété $K(\pi,1)$ pour $\ell.$
 \item Le morphisme $\phi_2:H^2(G_S^T(K)(\ell))\to H^2_{et}((X-S,T)_K)$ est surjectif et $cd\,G_S^T(K)(\ell)\leq 2.$
 \end{enumerate}

 Cette propriété est très forte, comme en témoigne les travaux de Schmidt à son propos, et elle permet de déduire certaines vertus pour les extensions maximales $S$-ramifiées, $T$-décomposées. 

 \subsubsection{Préparatifs}
 La première chose à faire est d'annuler le $\ell$-groupe des $T$-classes d'idéaux de $k.$
\begin{lemme}\label{lemme1} Soit $k$ un corps global et $T$ un ensemble non vide de places finies, telles que $\mathrm{pgcd}(\ell,\deg\tg,\ \tg\in T)=1$ dans le cas des corps de fonctions. Il existe $T_0$ tel que $_\ell Cl^{T\cup T_0}=0.$ De plus, on peut prendre $|T_0|=dim_\ell Cl^T\leq A_6\,g_k$ et $$\pi(T_0)\leq A_{9}g_k^2,$$ où  $A_{9}$ est une constante effective. 
\end{lemme}
\begin{preuve} Supposons que $g_k\geq 1,$ le cas $g_k=0$ étant clair. On procède par récurrence. On prend d'abord une place $\mathfrak{t}$ de $k$  telle que $\mathfrak{t}$ n'est pas totalement décomposée dans $K=k^{T,el}(\ell).$ D'après le théorème de densité de Cebotarev, on peut prendre une telle place $\mathfrak{t}$ telle que $\log\mathrm{N}\mathfrak{t}\leq c\log g_K.$ On recommence en remplaçant $T$ par $T\cup\{\mathfrak{t}\}$ et le lemme s'ensuit, en utilisant la majoration (utilisée précédemment) $dim_\ell Cl^T\leq A_{6} g_k$ puisque $\ell\geq 3$ (induisant $\log g_K\ll g_k$).
\end{preuve}

Nous aurons également besoin d'annuler les groupes de Kummer $V_Q^P:$ 
\begin{lemme} \label{lemme2}Soit $k$ un corps global. Soit $P$ un ensemble fini de places finies de $k.$ Alors il existe un ensemble fini de places finies $Q$ disjoint de $P,$ dont les places vérifient $\mathrm{N}\qg=1\ \mathrm{ mod }\ \ell$ pour tout $\qg\in Q,$ et 
$$V^P_{Q-\{\qg\}}(k)=0 \quad \forall\qg\in Q.$$ De plus on peut prendre $Q$ tel que:
$|Q|\leq 2(A_{6}g_k+\Phi_\mathbb{R}(k)+\Phi_\mathbb{C}(k)+|P|-1+\delta_\ell)$ et $$\pi(Q) \leq A_{10} |Q|\left( \log^{+}{|P|}+\log{\ell}+\log\{g_k+\delta_\Q n_k \log\ell+\pi(P)\}+\delta_P\,\ell\right),$$ où $A_{10}$ est effective.
\end{lemme}

\begin{preuve} Puisque $$V_Q^P=\ker\left\{V_\emptyset^P\to \prod_{\qg\in Q}k_\qg^\times/k_\qg^{\times \ell}\right\},$$ il s'agit de trouver deux places $\qg_{\alpha}$ telles que $\alpha\notin k_{\qg_{\alpha}}^{\times\ell}$ pour tout $\alpha$ dans une base de $V_\emptyset^P.$ D'après la théorie de Kummer, les restrictions des places de $k(\mu_\ell)$ qui ne sont pas totalement décomposées dans $k_\alpha=k(\mu_\ell,\sqrt[\ell]{\alpha})|k(\mu_\ell)$ conviennent. Comme il faut que deux telles places ne définissent pas la même place dans $k$ par restriction, on prend la seconde non conjuguée à la première par $Gal(k(\mu_\ell)|k).$ Il faut enfin les prendre hors de $P$ ce qui aura pour effet d'augmenter considérablement la borne.

Pour tout $\alpha\in V_\emptyset^P,$ $k_\alpha|k$ est par la théorie de Kummer non ramifiée hors de $P\cup\{\ell\},$ et de degré divisant $\ell(\ell-1)^{1-\delta_\ell}.$ Le genre de $k_\alpha$ est alors majoré par (de même que précédemment):
\begin{align*}
g(k_\alpha)&\leq \ell(\ell-1)^{1-\delta(k)}(g_k+\delta_\Q n_k \log\ell+\pi(P)).
\end{align*}
Appliquant deux fois le théorème de densité de Cebotarev en excluant les places de $P$ d'abord, puis $P$ et les conjugués de la place produite ensuite, on obtient deux places non conjuguées $\Qg_1,\Qg_2$ telles que $\alpha\notin k^{\times\ell}_{\qg_i}$ pour $i=1,2,$ où $\qg_i=\Qg_i\cap k.$ De plus $$\sum_{i=1,2}\log\mathrm{N}\Qg_i\ll \log{|P|}+\log{\ell}+\log\{g_k+\delta_\Q n_k \log\ell+\pi(P)\}+\delta_P\,\ell.$$ Notons qu'ici l'extension des constantes est au plus de degré $\ell,$ et, lorsque $_\ell Cl^P$ est trivial, il n'y en a pas, puisque $k_\alpha|k(\mu_\ell)$ est totalement ramifiée en une place de $P$ au moins.
On ajoute ainsi au plus $2\dim_\ell V^P_\emptyset(k)$ places. Comme $\dim_\ell V^P_\emptyset(k)\leq d_\ell(Cl^P)+\Phi_\mathbb{R}(k)+\Phi_\mathbb{C}(k)+|P|-1+\delta_\ell$ d'après les formules de Schafarevich, on obtient bien 
$$\pi(Q)\ll |Q|\left( \log^{+}{|P|}+\log{\ell}+\log\{g_k+\delta_\Q n_k \log\ell+\pi(P)\}+\delta_P\,\ell\right).$$
\end{preuve}

\subsubsection{A travers la preuve de \ref{Schmidt1}}
Ajoutons à présent une estimation quantitative au théorème suivant de Schmidt:
\begin{theo}[Schmidt] \label{schminter} Soit $T$ un ensemble fini de places de $k$ un corps global, non vide dans le cas des corps de fonctions. Soit $\ell\neq 2,$ $\ell\neq car(k).$ Alors il existe un ensemble fini $T_0$ de places de $k$ ainsi qu'un ensemble $S$ non vide dont les places $\pg$ vérifient $\N\pg=1\ \mod\ \ell$ telles que 
\begin{enumerate}
\item $S\cap(T\cup T_0)=\emptyset.$
\item $(X-S,T\cup T_0)$ vérifie la propriété $K(\pi,1)$ pour $\ell.$
\item Toute $\pg\in S$ se ramifie dans $k_S^{T\cup T_0}(\ell).$
\item $V_S^{T\cup T_0}=0.$
\end{enumerate}
\end{theo}

\begin{prop} On peut prendre dans le théorème \ref{schminter}, lorsque $\delta=0,$ et $pgcd(\ell,\deg\tg,\, \tg\in T)=1,$
\begin{enumerate}
\item $|T_0| \leq A_{6}g_k$ et $\pi(T_0) \leq A_{9}g_k^2$
\item $|S_0| \leq 2( h_k+\Phi_\mathbb{R}(k)+\Phi_\mathbb{C}(k)+|T'|-1+\delta_k)$ et 
\begin{align*}
\pi(S_0)& \leq A_{10}(A_{6}g_k+\Phi_\mathbb{R}(k)+\Phi_\mathbb{C}(k)+|T'|-1+\delta_k)\\
& \quad \times \left(\log^{+}|T\cup T_{0}|+\log\ell+\log\{g_k+\delta_\Q n_k \log\ell+\pi(T')\} \right)
\end{align*}
\item $|S| \leq 2|S_0|$ et
\begin{align*}
\pi(S)& \ll \pi(S_0)+|S_0|^2\log{\ell}\\
& \quad + |S_0|\left((h_k+|T'|+\Phi_\mathbb{R}(k)+\Phi_\mathbb{C}(k)-1)\log\ell+\pi'(T'\cup S_0)+\log^{+}{g_k}\right)
\end{align*}
\end{enumerate}
où $T'=T\cup T_0$ et $S_0$ est un sous-ensemble de $S$ construit de sorte que $$V^{T\cup T_0}_{S_0-\{\qg\}}(k)=0 \quad \forall\qg\in S_0.$$ 
\end{prop}
Dans l'application que nous avons en tête, $k$ sera égal à $\QQ,$ c'est pourquoi nous écrivons ce corollaire pour $\QQ.$
\begin{coro}
 Si $k=\QQ,$ $\delta=0,$ et $\mathrm{pgcd}(\ell,\deg\tg,\tg\in T)=1,$ on peut prendre $T_0$ et $S$ de sorte que :
\begin{align*}
(\QQ) \quad |T_0|&=\emptyset \\
|S|& \leq 4|T|\\
\pi(S)& \leq A_{11}\left(|T|\pi'(T)+|T|^2\log{\ell}\right),
\end{align*} pour une certaine constante effective $A_{11}.$
\end{coro}
  On peut obtenir par les mêmes estimations des résultats effectifs pour le cas $\delta=1,$ en regardant avec attention la preuve du théorème de Schmidt dans ce cas. Toutefois, on prend alors des places correspondant à chaque élément du groupe de Galois de $k_T^{el},$ c'est à dire que $|S|$ devient très grand, et il vaudra mieux jouer sur le nombre premier $\ell.$

\begin{preuve}
Tous les arguments algébriques de cette preuve sont dus à Schmidt, nous n'y ajoutons que les estimations des normes des places des ensembles intervenant. On fait la preuve de la proposition puis de son corollaire immédiat dans le même temps.

 Suivant Schmidt \cite{SGST}, on choisit d'abord $T_0$ tel que $T\cup T_0$  tel que $_\ell Cl^{T\cup T_0}=0.$ C'est possible d'après le lemme \ref{lemme1}. Prenons alors $S_0$ comme au lemme \ref{lemme2}, de sorte que $$V^{T\cup T_0}_{S_0-\{\qg\}}(k)=0 \quad \forall\qg\in S_0.$$ Posons $T'=T\cup T_0.$
On a alors: $|T'|\leq |T|+A_{6}g_k,$ 
$$\pi(T')\leq \pi(T)+A_{9}g_k^2,$$
$|S_0|\leq 2( h_k+\Phi_\mathbb{R}(k)+\Phi_\mathbb{C}(k)+|T'|-1+\delta_k), $ et 
\begin{align*}
\pi(S_0)\leq & A_{10}(A_{6}g_k+\Phi_\mathbb{R}(k)+\Phi_\mathbb{C}(k)+|T'|-1+\delta_k)\\ 
& \times \left(\log^{+}|T\cup T_{0}|
 +\log\ell+\log\{g_k+\delta_\Q n_k \log\ell+\pi(T')\}\right)
\end{align*}
Dans le cas où $k=\QQ,$ on a déjà $_\ell Cl^{T}=0,$ ainsi on prendra $T_0=\emptyset$ dans ce cas. Alors 
$ (\QQ)\quad |S_0|\leq 2 |T|,$ et 
$$(\QQ)\quad\pi(S_0)\leq 2\,A_{10}|T|\left(\log^{+}|T|+\log\ell+\log\{\log\ell+\pi(T)\} \right)$$ et donc 
$$(\QQ)\quad\pi(S_0)\ll |T|\left(\log^{+}{\pi(T)}+\log\ell\right)$$

On écrit $S_0=\{\pg_1,\dots,\pg_m\}.$ Du fait de $_\ell Cl^{T'}=0,$ l'application $(\mathrm{v}_\pg)_\pg:k^\times\to  \bigoplus_{q\notin T'}\mathbb{Z}$ induit alors la suite exacte $$0\to E_{k,T'}/\ell\to k^\times/k^{\times \ell}\to \bigoplus_{q\notin T'}\mathbb{Z}/\ell\mathbb{Z}\to 0.$$ Pour une place $\pg\notin T',$ on peut alors considérer un élément $s_\pg\in k^\times/k^{\times \ell}$ défini par $\mathrm{v}_\pg(s_\pg)=1 \mod \ell$ et $\mathrm{v}_\qg(s_\pg)=0 \mod \ell$ pour tout $\qg\notin T'\cup\{\pg\}.$ Il est défini à $E_{k,T'}/\ell$ près. Pour tout $i=1\dots m,$ on considère $s_{\pg_i}$ correspondant à $\pg_i$  et on le notera plus simplement $s_i.$

Rappelons de plus deux lemmes dus à Schmidt \cite{SGST}:
\begin{lemme}[Schmidt] \label{lemme52} Soit $T$ un ensemble fini de places de $k,$ non vide dans le cas des corps de fonctions, tel que $_\ell Cl^T(k)=0.$ Soit $\qg\notin T\cup P_\ell$ une place totalement décomposée dans $k(\sqrt[\ell]{E_{k,T}})|k.$ Alors l'extension $k^{T,el}_{\{\qg\}}|k$ est cyclique d'ordre $\ell$ et $\qg$ se ramifie dans cette extension. Enfin un idéal premier $\pg\notin T\cup\{\qg\}$ de norme $\N\pg=1\mod\ell$ se décompose dans  $k^{T,el}_{\{\qg\}}$ si et seulement si $\qg$ se décompose totalement dans $k(\sqrt[\ell]{E_{k,T}},\sqrt[\ell]{s_\qg})|k(\sqrt[\ell]{E_{k,T}}).$
\end{lemme}

\begin{lemme}[Schmidt]\label{lemme64}
Si $\delta=0$ et $S=\{\pg_1,...,\pg_n\}$ est un ensemble fini de places finies avec $\N(\pg_i)=1 \mod{\ell}.$ Posons $s_i=s_{\pg_i}.$ Alors les extensions $k(\mu_\ell,\sqrt[\ell]{s_1},\dots,\sqrt[\ell]{s_n})$ et $k_S^{T,el}(\mu_\ell)$ sont linéairement disjointes sur $k(\mu_\ell).$
\end{lemme}
Remarquons également qu'aucune des extensions qu'on considère ici n'admet d'extension des constantes au dessus de $k(\mu_\ell),$ d'après les conditions sur les places de $T,$ l'hypothèse $\delta_\ell=0$ et les propriétés des extensions de Kummer (chacune sera ramifiée ou $T$ décomposée).

On va maintenant ajouter des places à $S_0,$ choisies de la manière suivante. Soient $\Pg_1,\dots,\Pg_m$ des prolongements de $\pg_1,\dots,\pg_m$ à $k(\mu_\ell).$ On considère, pour une place $\Qg$ de $k(\mu_\ell),$ et $a\in\{1,\dots,m\}$ la propriété $(B_a):$
\begin{enumerate}
\item $\Qg\notin T'(k(\mu_\ell)),$
\item $Frob_\Qg\notin \mathcal{I}_{\Pg_a}\subset G(k_{S_0}^{T',el}(\mu_\ell)|k(\mu_\ell)),$
\item Pour tout $b\neq a,$ $\Qg$ se décompose dans $k(\mu_\ell,\sqrt[\ell]{s_b})|k(\mu_\ell),$
\item $\Qg$ est inerte dans $k(\mu_\ell,\sqrt[\ell]{s_a})|k(\mu_\ell),$
\item $\Qg$ est totalement décomposée dans $k(\sqrt[\ell]{E_{k,T'}})|k(\mu_\ell).$
\end{enumerate}

Cette propriété est ainsi indépendante du choix des $s_i.$ On construit $\Qg_1,...,\Qg_m$ de la façon suivante: on prend $\Qg_1\in P(k(\mu_\ell))$ vérifiant $B_1.$ On pose $\qg_1=\Qg_1\cap k.$ On voit alors que $k_{\qg_1}^{T',el}$ est cyclique d'ordre $\ell$ et ramifiée en $\qg_1.$ On choisit alors par récurrence les places $\Qg_2,...,\Qg_m$ (et on pose $\qg_a=\Qg_a\cap k$) de sorte que
\begin{enumerate}
 \item $\Qg_a$ vérifie la propriété $(B_a)$ 
 \item $\Qg_a$ est, pour $b<a,$ décomposée dans $k_{\qg_1}^{T',el}(\mu_\ell)|k(\mu_\ell)$ et $k(\mu_\ell,\sqrt[\ell]{s_{\qg_b}}).$
\end{enumerate}

Un tel choix est possible du fait du lemme \ref{lemme64}. Estimons alors la taille des $\Qg_a.$ Choisir une telle place revient à imposer des condition sur son Frobenius dans le compositum des extensions linéairement indépendantes (d'après les lemmes \ref{lemme52} \ref{lemme64}) $k(\mu_\ell,\sqrt[\ell]{s_1},\dots,\sqrt[\ell]{s_n},\sqrt[\ell]{s_{\qg_1}},\dots,\sqrt[\ell]{s_{\qg_{a-1}}}),$ $k(\sqrt[\ell]{E_{k,T'}}),$ $k_{S_0}^{T',el}(\mu_\ell)$ et les $k^{T',el}_{\{q_1\}},\dots,$ $k^{T',el}_{\{q_{a-1}\}}.$ 
Estimons alors le genre de ce compositum $L_a.$ Comme toutes les extensions (à partir de $k(\mu_\ell)$) sont des $\ell$-extensions, on calculera plutôt le logarithme en base $\ell$ des degrés (nommé $\ell$-degré). Si $K/k$ est une $\ell$-extension, on posera $[K:k]_\ell=\log_\ell [K:k].$

Commençons par calculer $[k_{S_0}^{{T'},el}(\mu_\ell):k(\mu_\ell)]_\ell.$ 
Comme $V_{S_0}^{T'}=0,$ on a $$h^1(k_S^{T'}(\ell))=1+|S_0|-\Phi_\mathbb{R}(k)-\Phi_\mathbb{C}(k)-|T'|$$ (voir \cite[10.7.10]{NCG2}), et donc 
$$[k_{S_0}^{{T'},el}(\mu_\ell):k(\mu_\ell)]_\ell\leq 2h_k+|T'|+\Phi_\mathbb{R}(k)+\Phi_\mathbb{C}(k)-1.$$ 


$k(\sqrt[\ell]{E_{k,T'}})|k(\mu_\ell)$ étant une extension de $\ell$-degré $\dim E_{k,T'}/\ell= r-1+|T'|,$ on en déduit que 
\begin{align*}
[L_a:k]|&(\ell-1)\ell^{[k_{S_0}^{{T'},el}:k]_\ell+[k(\sqrt[\ell]{E_{k,T'}}):k(\mu_\ell)]_\ell+|S_0|+2(a-1)}.
\end{align*}
D'après la théorie de Kummer, l'extension $L_a$ est non ramifiée hors de $R=S_0\cup T'\cup \{\qg_1,\dots,\qg_{a-1}\}\cup\{\ell\}.$ De plus, seules les places au-dessus de $\ell$ dans $k$ peuvent être sauvagement ramifiées dans $L_a/k,$ on obtient donc, d'après le lemme \ref{genre} :
\begin{align*}
g^\ast_{L_a}&\leq [L_a:k] \left(g^\ast_k+\frac{1}{2}\sum_{v\in R}\log{\mathrm{N}v}+\frac{\delta_\Q n_k}{2}[L_a:k(\mu_\ell)]_\ell\log{\ell}\right):=\bar{g}_{L_a} \\
\end{align*}

D'après le théorème de densité de Cebotarev, on peut trouver $\Qg_a$ comme on le demande, avec en plus $\log \N\Qg_a\ll \log^{+}|T'|+\log (n_{L_a}\bar{g}_{L_a}).$ 
On a alors 

\begin{align*}
\log \N\Qg_a &\ll \log^{+}|T'|+([k_{S_0}^{{T'},el}(\mu_\ell):k(\mu_\ell)]_\ell+[k(\sqrt[\ell]{E_{k,T'}}):k(\mu_\ell)]_\ell\\
& +|S_0|+2(a-1)+1)\log\ell +\log \Bigg\{g^\ast_k+\frac{1}{2}(\pi(T'\cup S_0)+\log\ell+\sum_{i<a}\log\N\Qg_a\\
&  +\delta_\Q n_k\left([k_{S_0}^{{T'},el}(\mu_\ell):k(\mu_\ell)]_\ell+[k(\sqrt[\ell]{E_{k,T'}}):k(\mu_\ell)]_\ell+|S_0|+2(a-1)\right)\log\ell\Bigg\}\\
&\ll \log^{+}{|T'|}+(h_k+|T'|+\Phi_\mathbb{R}(k)+\Phi_\mathbb{C}(k)-1+|S_0|+(a-1))\log\ell \\
& \quad +\pi'(T'\cup S_0)+\log^{+}{g_k}+\log\sum_{i<a}\log\N\Qg_a
\end{align*}

Posons $A= \log^{+}{|T'|}+(h_k+|T'|+\Phi_\mathbb{R}(k)+\Phi_\mathbb{C}(k)-1+|S_0|)\log\ell+\pi'(T'\cup S_0)+\log^{+}{g_k},$ et $X_0=0,$ $X_k=\sum_{i\leq k}\log\N\Qg_k$ pour $k\geq 1.$ On voit que $A\geq 2\log 3.$
On a alors, pour $k\geq 1,$ $$X_k\ll A k+k^2\log\ell+k\log X_{k-1}. $$ Montrons alors que  $X_k\ll Ak+k^2\log\ell.$ 

Partons de l'inégalité $X_k\leq c(A k+k^2\log \ell+k\log X_{k-1})$ pour une certaine constante effective $c\geq 2.$ Considérons $Y_k=X_k/c.$ On a alors $Y_k\leq (A+\log{c})k+k^2\log\ell+k\log Y_k.$ Posons $B=A+\log{c}.$

Montrons par récurrence que $$Y_k\leq 8\left(B k+k^2\log\ell\right) .$$ En effet, pour $k=0$ c'est clair. Supposons la vérifiée pour $k-1,$ on a alors:
\begin{align*}
Y_k &\leq Bk+k\log{B}+k^2\log\ell+3 k\log{k}+k(\log\log\ell+ 4\log{2})\\
& \leq 2Bk+8 k^2\log\ell. \\
\end{align*}

On en déduit donc que: 
\begin{align*}
X_m& \ll \left( \log^{+}{|T'|}+(h_k+|T'|+\Phi_\mathbb{R}(k)+\Phi_\mathbb{C}(k)-1)\log\ell+\pi'(T'\cup S_0)+\log^{+}{g_k}\right)|S_0|\\
& \quad +|S_0|^2\log{\ell}.
\end{align*} 

On pose alors $S=S_0\cup\{\qg_1,\dots,\qg_m\}.$ On obtient que $|S|=2|S_0|$ et $\pi(S)\leq \pi(S_0)+X_m.$ On obtient alors la preuve de la proposition.

Dans le cas de $\QQ,$ on a alors:
$$(\QQ)\quad |S|\leq 4|T|,$$ et 

\begin{align*}
(\QQ) \quad \pi(S) & \ll |T|(\pi'(T)+\log{\ell})+2|T|\left(\log^{+}|T|+|T|\log\ell +\pi'(T)+\pi'(S_0)\right))\\
&\quad +4|T|^2\log\ell \\
& \ll |T|\pi'(T)+|T|^2\log{\ell}
\end{align*}
On peut alors montrer que cet ensemble $S$ convient (voir la preuve du théorème \cite[6.1]{SGST}). 
\end{preuve}

\subsubsection{Version quantitative du théorème \ref{Schmidt1}.} 

\begin{prop} \label{Schmidt1exp} Dans le cas de $\QQ$ et lorsque $\delta=0,$ on peut prendre dans le théorème \ref{Schmidt1}
\begin{align*}
(\QQ)\quad S_0 &\leq 4|S\cup T|\\
\pi(S_0)& \leq  A_{11}\left(|T\cup S|\pi'(T\cup S)+|T\cup S|^2\log{\ell}\right)
\end{align*}
\end{prop}

\begin{preuve} En effet, d'après \cite[7.1]{SGST}, il suffit de trouver $S_0$ telle que $(X-S_0,S\cup T)$ ait la propriété $K(\pi,1)$ pour $\ell.$ On obtient donc toutes les informations voulues de celles du paragraphe précédent, avec $T$ remplacé par $T\cup S.$ 
\end{preuve}

\subsection{Application aux corps globaux infinis}

Nous allons à présent tirer de ces résultats des informations sur les invariants de corps globaux infinis.

\begin{theo} \label{QST} Soit $\ell$ un nombre premier impair ne divisant pas $r-1$ dans le cas des corps de fonctions. Soit $T$ et $I$ deux ensembles finis disjoints de places finies de $\QQ.$ Alors il existe un ensemble fini $S$ de places de $\QQ$ de norme congrue à $1\mod{\ell}$ tel que $\QQ_{S}^T(\ell)$ a les propriétés suivantes:
\begin{enumerate}
\item Pour tout $\pg\in Pl_f(\QQ),$ $\phi_{\pg,\N\pg^m}=0$ si $m\geq 2,$
\item Pour tout $\pg\in T,$ $\phi_{\pg,\N\pg}=\phi_\infty>0,$
\item Pour tout  $\pg\in I\cup S,$ $\phi_{\pg,\N\pg}=0.$
\item  Dans le cas des corps de nombres, \begin{displaymath} \frac{2}{\pi(S)}\leq \phi_\mathbb{R}\leq \frac{4}{\pi(S)}\end{displaymath} et $\phi_\mathbb{R}=\phi_\infty.$
\item On peut prendre $S$ tel que $|S|\leq 4|T+1|+1$ et $$\pi(S) \ll  |T|(\pi'(T)+\log^{+}|I|+\log^{+}\pi'(I))+\pi'(I)+(|T|^2+|I|+1)\log{\ell}.$$
\end{enumerate}

\end{theo}

\begin{preuve} On va prendre $S$ tel que la $\ell$-dimension cohomologique de $Gal(\QQ_{S}^T(\ell)/\QQ)$  soit finie. Ainsi ce groupe n'aura pas de torsion, et donc pas d'invariants non nuls hormis ceux de degré $1$ (correspondant à $m=1$). De plus $T$ sera totalement décomposé, et on aura donc  $\phi_{\pg,\N\pg}=\phi_\infty>0,$ pour tout $\pg\in T.$  $\phi_\infty$ est strictement positif car l'extension est modérément ramifiée, non ramifiée hors d'un ensemble fini de places. Enfin il faudra prendre $S$  de sorte que les places de $I$ ne soient pas totalement décomposées. On peut le faire et même estimer la norme des places ainsi ajoutées, au moyen du théorème de Grunwald-Wang et des résultats de la section précédente. Enfin l'estimation concernant $\phi_\infty$ provient de la majoration suivante, valable pour toute extension galoisienne (voir \cite{Lth}). Soit $\{K_i\}_{i\in \mathbb{N}}$ une tour représentant $\KK=\QQ_{S}^T(\ell),$ que l'on suppose construite. Toutes les places de $S$ sont alors ramifiées dans $\KK$ puisqu'elle réalise même l'extension maximale locale en les places de $S.$ On a donc, pour $i$ suffisamment grand pour que toutes les places de $S$ soient ramifiées dans $K_i/\QQ,$
$$ g^\ast_{K_{i}}= g^\ast_\QQ n_{K_i}+\frac{1}{2}\sum_{\pg\in S}(e_\pg-1)g_\pg f_\pg\log\N\pg,$$
et donc, puisque $e_\pg\geq 2,$ $e_\pg/2\leq e_\pg-1\leq e_\pg,$ d'où l'on déduit:
$$\frac{1}{4}\pi(S)-\theta \leq \frac{g^\ast_{K_{i}}}{n_{K_i}}\leq \frac{1}{2}\pi(S)-\theta,$$ où $\theta$ vaut 1 dans le cas des corps de fonctions, $0$ sinon, et donc, dans le cas des corps de nombres, on obtient:  
$$\frac{2}{\pi(S)}\leq \phi_\infty\leq \frac{4}{\pi(S)}.$$

Estimons alors $S.$ Soit $s$ une place telle qu'il existe une extension galoisienne de $\QQ$ d'exposant $\ell,$ non ramifiée hors de $s$ où toutes les places de $I$ sont inertes, et où les places de $T$ sont totalement décomposées. Une telle extension existe d'après le théorème de Grunwald-Wang, et on peut prendre $s$ telle que
$\log\N s\ll \log \bar{g}_L$ (avec les notations du $2.1$) où
$$\bar{g}_L=(\ell-1) \ell^{\#T+\#I} \left(\frac{1}{2}\pi(T\cup I)+\delta_\Q\frac{\log{\ell}}{2} (\#I+\#T+ 1)\right).$$ On a alors:
$$\log\N s \ll (1+|T\cup I|)\log{\ell} +\pi'(T\cup I).$$
On considère alors $S_0$ tel que $cd_\ell \KK=2$ avec $\KK=\QQ_{S_0\cup\{s\}}^T(\ell).$ D'après la proposition \ref{Schmidt1exp}, on peut prendre $S_0$ tel que:
$|S_0|\leq 4(|T|+1)$ et 
$$\pi(S_0) \ll  |T+1|(\pi'(T)+\log\log\N s)+(|T|+1)^2\log{\ell}.$$ On en déduit que 
$$\pi(S_0) \ll  |T+1|(\pi'(T)+\log^{+}|I|+\log^{+}\pi'(I))+|T+1|^2\log{\ell}.$$ On pose $S=S_0\cup\{s\},$ et on obtient le résultat escompté:
$$\pi(S) \ll  |T+1|(\pi'(T)+\log^{+}|I|+\log^{+}\pi'(I))+\pi'(I)+(|T|^2+|I|+1)\log{\ell}.$$
\end{preuve}

\subsection{Preuve de la proposition \ref{propB}}
Notons $\KK=\QQ_S^T(\ell).$ Prouvons le résultat dans le cas des corps de nombres et sous $(GRH),$ les résultats sans cette hypothèse et pour les corps de fonctions se déduisant immédiatement en utilisant l'inégalité correspondante. On utilisera cependant les $\phi_{\pg,q}$ plutôt que les $\phi_q,$ alourdissant ainsi les notations, de sorte que la preuve pour les corps de fonctions soit exactement la même. On a démontré précédemment que $2/\pi(S)\leq\phi_\infty(\KK)\leq 4/\pi(S).$ De plus, pour les places totalement décomposées, $\phi_{\pg,\N\pg}=\phi_\infty.$ $\KK$ est totalement réel dans le cas des corps de nombres; on a alors, d'après les inégalités fondamentales de Tsfasman-Vl\u adu\c t, 
$$\sum_{\pg\in D}\frac{\phi_{\pg,\N\pg}\log{\N\pg}}{\sqrt{\N\pg}-1}\leq 1-(\log{\sqrt{8\pi}}+\frac{\pi}{4}+\frac{\gamma}{2})\phi_\infty,$$ et donc
\begin{align*}
\sum_{\pg\in D}\frac{\log{\N\pg}}{\sqrt{\N\pg}-1}& \leq \frac{1}{\phi_\infty}-(\log{\sqrt{8\pi}}+\frac{\pi}{4}+\frac{\gamma}{2}),\\
& \leq \frac{\pi(S)}{2}-\log{\sqrt{8\pi}}-\frac{\pi}{4}-\frac{\gamma}{2}
\end{align*}

\section{Preuve du théorème \ref{princ}}

On suppose encore ici \emph{GRH} pour ce qui est des majorations. Soient $P=\{\pg_1,\dots,\pg_n\}$ un ensemble de $n$ places finies de $\QQ,$ $d_{i,1},\dots,d_{i,n_i}$ $n_i$ entiers naturels donnés pour tout $i=1\dots n.$ Soit $I$ un ensemble de places finies disjoint de $P.$

D'après les résultats du second paragraphe, il existe une extension de $\QQ$ telle que:
  \begin{enumerate}
\item $\Phi_{\pg_i,\N\pg^{d_{i,j}}}(L)=\frac{[L:\QQ]}{n_i\,d_{i,j}}.$
\item Pour tout $\Pg\in P_L$  au-dessus d'un des $\pg_i,$ il existe $j$ tel que $\mathrm{N}\pg=\N\pg_i^{d_{i,j}}.$
\item Il existe une fonction $f$ de $P$ et $N=\mathrm{ppcm}(n_i)_i\mathrm{ppcm}(d_{i,j})_{i,j}$ telle que $$g_L/[L:\QQ]\leq f(P,N).$$
De plus, sous l'hypothèse de Riemann généralisée et si $r$ est premier avec $N,$ $f(P,N)$ peut être prise ainsi:
\begin{align*}
(GRH)\quad f(P,N) &\ll A_7^{\Omega(N)}\left\{(1+|P|)\log{N} +\pi'(P)+\delta_\F N^2\right\}, \end{align*}
   \end{enumerate}
 Choisissons $\ell$ à présent. Dans le cas des corps de nombres, on prend $\ell=3.$ Dans le cas des corps de fonctions, le choix doit être fait de sorte à proscrire les extensions des constantes. Il ne doit donc pas diviser $r-1$ ni $a_P=\mathrm{pgcd}(\deg\pg,\ \pg\in P).$ On prend alors $\ell$ vérifiant $r+a_P\leq \ell\leq 2(r+a_P).$

 Il existe $S$ tel que $\KK=\QQ_{S}^P(\ell)$ vérifie:  
\begin{enumerate}
\item pour tout $\pg\in Pl_f(\QQ),$ $\phi_{\pg,\N\pg^m}=0$ si $m\geq 2,$
\item pour tout $\pg\in P,$ $\phi_{\pg,\N\pg}=\phi_\infty>0,$
\item pour tout  $\pg\in I\cup S,$ $\phi_{\pg,\N\pg}=0,$
\item dans le cas des corps de nombres \begin{displaymath} \frac{2}{\pi(S)}\leq \phi_\mathbb{R}=\phi_\infty\leq \frac{4}{\pi(S)}.\end{displaymath} 
\item  $|S|\leq 4|P|+5$ et $\frac{1}{2}\pi(S)\leq g(P,I,\ell),$ avec $$g(P,I,\ell) \ll  |P|(\pi'(P)+\log^{+}|I|+\log^{+}\pi'(I))+\pi'(I)+(|P|^2+|I|+1)\log{\ell}.$$
\end{enumerate}

On considère alors le compositum $L\KK.$  L'extension $L\KK/L$ est non ramifiée hors de $S_L,$ modérément ramifiée. Ainsi $\phi_\infty(L\KK)>0.$ De plus, toutes les places de $P_L$ y sont totalement décomposées. Pour tout $\pg\in S\cup I,$ on a bien $\phi_{\pg,q}(L\KK)=0,$ pour tout $q,$ puisque $$\sum_q \phi_{\pg,q}(L\KK)\log_{\N\pg}{q} \leq \sum_q \phi_{\pg,q}(\KK)\log_{\N\pg}{q}$$ (car $\KK\subset L\KK$).

On a, pour les places $\pg_i\in P$ et pour une tour $K_i$ représentant $\KK,$  
$$\text{pour tout } j=1\dots n_i,\quad\Phi_{\pg_i,\N\pg_i^{d_{i,j}}}(LK_i)=\frac{[LK_i:\QQ]}{n_i\,d_{i,j}},$$ d'où l'on déduit que 
$$\text{pour tout } j=1\dots n_i,\quad\phi_{\pg_i,\N\pg^{d_{i,j}}}(L\KK)=\frac{\phi_\infty(L.\KK)}{n_i\,d_{i,j}}>0.$$
On voit alors que $\delta(L\KK)\leq 1-\varepsilon,$ où
$$\varepsilon=\phi_\infty(L.\KK)\sum_{i=1}^{n}\frac{\log{\N\pg_i}}{n_i}\sum_{j=1}^{n_i}\frac{1}{\N\pg_i^{\frac{d_{i,j}}{2}}-1}+\delta_\Q\left(\log{8\pi}+\frac{\pi}{4}+\frac{\gamma}{2}\right).$$
Reste donc à minorer $\phi_\infty(L\KK),$ majorons alors le quotient $g_{LK_i}/n_{LK_i}.$
$LK_i/L$ étant modérément ramifiée, non ramifiée hors de $S_L$ on a donc:
$$g^\ast_{LK_i}\leq[L:K_i]g^\ast_L+\frac{1}{2}[LK_i:L]\sum_{\Pg\in S_L}\log\N\Pg.$$
Comme $$\sum_{\Pg\in S_L}\log\N\pg=\sum_{\pg\in S}\log\N\pg\sum_{\Pg\in S_L}f_{\Pg|\pg}\leq [L:\QQ]\pi(S),$$ on obtient alors 
$$\frac{g^\ast_{LK_i}}{[LK_i:\QQ]}\leq \frac{g^\ast_L}{[L:\QQ]}+\frac{1}{2}\pi(S).$$
On en déduit que $\phi_\infty\geq (g(P,I,\ell)+f(P,N))^{-1},$ ce qui termine la preuve.
\section{Deux cas particuliers}
On se propose à présent de donner les estimations dans les deux cas particuliers les plus significatifs.

\subsection{Cas particulier $N=1$}

Dans ce dernier paragraphe, nous allons estimer le défaut qu'on obtient dans le cas particulier où on prend pour $T=\{2, 3,\dots,p_n\}$ les $n$ plus petits nombres premiers. On peut alors construire $S$ tel que $\Q_S^T(3)$ ait une dimension cohomologique $2$ et vérifiant les conditions du théorème \ref{QST}. Considérons alors le défaut $\delta$ de $\Q_S^T(3)$ sous $GRH.$ Notons $D$ l'ensemble des places de $\Q$ totalement décomposées dans $\Q_S^T(3).$ Comme seules les places totalement décomposées $\pg$ ont au moins un $\phi_{\pg,q}>0,$ on a: 

On a $$\delta=1-\sum_{\pg\in D}\frac{\phi_{\pg,\N\pg}\log{\N\pg}}{\sqrt{\N\pg}-1}-(\log{\sqrt{8\pi}}+\frac{\pi}{4}+\frac{\gamma}{2})\phi_\infty,$$ et donc 
\begin{align*}
\delta &=1-\phi_\infty\left(\sum_{\pg\in D}\frac{\log{\N\pg}}{\sqrt{\N\pg}-1}+\log{\sqrt{8\pi}}+\frac{\pi}{4}+\frac{\gamma}{2}\right),\\
\delta &\leq 1-\frac{2}{\pi(S)}\left(\sum_{p\leq p_n}\frac{\log{p}}{\sqrt{p}-1}+\log{\sqrt{8\pi}}+\frac{\pi}{4}+\frac{\gamma}{2}\right)
\end{align*}

On a :
$\pi(S)\ll |T|\pi'(T)+(|T|^2+1)\log{3}.$
On est alors amené à estimer $\log\sum_{p\leq p_n}\log{p}.$ La somme $\vartheta(x)=\sum_{p\leq x}\log{p}$ est bien connue et vaut $\vartheta(x)=x+o(x).$ Comme $p_n=n\log{n}(1+o(1))$ on en déduit que $\pi'(T)=\log{n}+\log\log{n}+o(1).$ Ainsi $\pi(S)\ll n^2.$
Reste à calculer $S_n=\sum_{p\leq p_n}\frac{\log{p}}{\sqrt{p}-1}.$ Intégrant par partie on obtient: 
$$S_n\geq \frac{\vartheta(p_n)}{\sqrt{p_n}}+\int_{t=2}^{p_n} \frac{\vartheta(t)dt}{2t\sqrt{t}}.$$
Comme $\vartheta(x)\gg x,$  on a:
$$S_n\gg \sqrt{p_n}\gg c\sqrt{n\log{n}},$$ puisque $p_n\gg n\log{n}.$ 
On en déduit qu'il existe une constante effective $c,$ telle que $\delta\leq 1-c\,n^{-3/2}\sqrt{\log{n}}.$ 

Remarquons que ce résultat est bien plus faible que celui obtenu dans \cite{Lth}, utilisant des extensions de degré $2.$

\subsection{Cas Particulier d'une place $p$ et de $n$ degrés}

Dans ce paragraphe on va considérer un premier $p$ et $d_1=1,\dots, d_n=n.$ Alors il existe un corps global infini $\mathcal{K}/\Q$ tel que, pour tout $k=1\dots n,$ $$\phi_{p^k}=\frac{\phi_\mathbb{R}}{kn}>0.$$ De plus, son défaut vérifie:
$$(GRH)\quad \delta\leq 1-\frac{h(p,n)}{f(p,n)+g(p)},$$ où
$$h(p,n)\geq \frac{\log{p}}{n\sqrt{p}}+\log{8\pi}+\frac{\pi}{4}+\frac{\gamma}{2},$$
$$f(p,n)\ll a^{\Omega(n\,\mathrm{ppcm}(k)_{k=1\dots n})}(n\log{n}+\log\log{p}), \text{ et } g(p)\ll \log\log{p}.$$

Reste alors à évaluer $\Omega(\mathrm{ppcm}(k)_{k=1\dots n})=\sum_{p\leq n}m_p$ où la somme est prise sur les nombre premiers, et $m_p$ est le plus grand entier tel que $p^{m_p}\leq n.$ Comme 
$$\sum_{p\leq n}m_p\leq \log{n}\sum_{p\leq n}\frac{1}{\log{p}}\ll \log{n} \frac{n}{\log^2{n}},$$ on voit que:
$$\delta\leq 1-\frac{c}{A^{n/\log{n}}(n\log{n}+\log\log{p})},$$ pour $c,A$ deux constantes effectives.

\bibliographystyle{alpha-fr}
\bibliography{biblio}

\end{document}